\documentclass[11pt]{amsart}

\newif\ifdraft\draftfalse

\usepackage{amssymb, enumerate, xspace, graphicx, url}
\usepackage[mathscr]{eucal}
\usepackage[all]{xy}
\SelectTips{cm}{}

\ifdraft\usepackage[notcite, notref]{showkeys}\fi

1

\numberwithin{subsection}{section}

\allowdisplaybreaks[1]

\newenvironment{enumeratea}
{\begin{enumerate}[\upshape (a)]}
{\end{enumerate}}

\newtheorem*{namedtheorem}{\theoremname}
\newcommand{\theoremname}{testing}

\newcommand\Fields{\operatorname{Fields}}

\newcommand\zp{\operatorname{Points}}
\newcommand\points{\zp}
\newcommand\colim{\operatorname{colim}}
\newcommand\Sets{\operatorname{Sets}}
\newcommand\sets{\operatorname{Sets}}

\newcommand\Aff{\operatorname{Aff}}
\newcommand\Alg{\operatorname{Alg}}

\newcommand\chr{\operatorname{char}}
\newcommand\Ind{\operatorname{Ind}}

\newcommand\Orb{\operatorname{\bf Orb}}
\newcommand\Spec{\operatorname{Spec}}
\newcommand\Pf{\operatorname{Pf}}

\newtheorem{theorem}[equation]{Theorem}
\newtheorem{proposition}[equation]{Proposition}
\newtheorem{proposition-definition}[equation]{Proposition-Definition}
\newtheorem{corollary}[equation]{Corollary}
\newtheorem{lemma}[equation]{Lemma}
\newtheorem{conjecture}[equation]{Conjecture}

\theoremstyle{definition}
\newtheorem{definition}[equation]{Definition}

\newtheorem{example}[equation]{Example}
\newtheorem{examples}[equation]{Examples}
\newtheorem{remark}[equation]{Remark}
\newtheorem{question}[equation]{Question}

\newtheorem{para}[equation]{\relax}

\theoremstyle{remark}
\newtheorem*{claim}{Claim}
\newtheorem{Claim}[equation]{Claim}

\numberwithin{equation}{section}

 \newcommand\cB{\mathcal{B}}
 \newcommand\cD{\mathcal{D}}
 \newcommand\cF{\mathcal{F}}
 \newcommand\cH{\mathcal{H}}
\newcommand\cI{\mathcal{I}} 
 
\newcommand\cM{\mathcal{M}} 
\newcommand\cO{\mathcal{O}}

\newcommand\cU{\mathcal{U}} 
 \newcommand\cX{\mathcal{X}}
\newcommand\cY{\mathcal{Y}} 

\renewcommand\AA{\mathbb{A}} 
\newcommand\CC{\mathbb{C}} 
 
\newcommand\GG{\mathbb{G}}

 \newcommand\PP{\mathbb{P}}
\newcommand\QQ{\mathbb{Q}}

 \newcommand\ZZ{\mathbb{Z}}

 \newcommand\bH{\mathbf{H}}

 \newcommand\bX{\mathbf{X}}

\newcommand\rC{\mathrm{C}}

\newcommand\rM{\mathrm{M}} 
\newcommand\rO{\mathrm{O}} 
 
 \newcommand\rT{\mathrm{T}}
 
\newcommand\rW{\mathrm{W}}

\newcommand\rma{\mathrm{a}} 
 
\newcommand\rme{\mathrm{e}} 
 \newcommand\rmh{\mathrm{h}}

\newcommand\rmm{\mathrm{m}}

\newcommand\fM{\mathfrak{M}}

 \newcommand\frp{\mathfrak{p}}

\newcommand\arr{\ifinner \to\else\longrightarrow\fi}

\newcommand\arrto{\ifinner\mapsto\else\longmapsto\fi}
\newcommand\larr{\longrightarrow}

\newcommand{\xarr}{\xrightarrow}
\newcommand{\darr}{\dashrightarrow}

\renewcommand\H{\operatorname{H}}

\newcommand\op{^{\mathrm{op}}}
\newcommand{\dual}{^{\vee}}
\newcommand{\inv}[1]{{#1}^{*}/{#1}^{*n}}

\newcommand\eqdef{\overset{\mathrm{\scriptscriptstyle def}} =}

\newcommand\into{\hookrightarrow}

\def\displaytimes_#1{\mathrel{\mathop{\times}\limits_{#1}}}

\def\displayotimes_#1{\mathrel{\mathop{\bigotimes}\limits_{#1}}}

\newcommand\Mor{\operatorname{Mor}}
\newcommand\Br{\operatorname{Br}}

\newcommand\ind{\operatorname{ind}}
\newcommand\rank{\operatorname{rank}}
\newcommand\red{\operatorname{red}}

\newcommand\ord{\operatorname{ord}}

\newcommand\spec{\operatorname{Spec}}

\newcommand\generate[1]{\langle #1 \rangle}
\newcommand\pfister[1]{{\ll #1 \gg}}

\newcommand\id{\mathrm{id}}

\newcommand\pr{\operatorname{pr}}

\newdir{ >}{{}*!/-5pt/@{>}}

\newcommand\double{\mathbin{\rightrightarrows}}

\newcommand\doublelong[2]{\mathbin{\xymatrix{{}\ar@<3pt>[r]^{#1}
\ar@<-3pt>[r]_{#2}&}}}

\newcommand{\underaut}
{\mathop{\underline{\mathrm{Aut}}}\nolimits}

\newlength{\ignora}
\newcommand{\hsmash}[1]{\settowidth{\ignora}{#1}#1\hspace{-\ignora}}

\newcommand{\ed}{\operatorname{ed}}

\newcommand{\dm}{Deligne--Mumford\xspace}

\newcounter{steps}

\newcommand{\trdeg}{\operatorname{tr\,deg}}
\newcommand{\mmu}{\boldsymbol{\mu}}

\newcommand{\GL}{\mathrm{GL}}
\newcommand{\SL}{\mathrm{SL}}
\newcommand{\PGL}{\mathrm{PGL}}

\newcommand\Spin{\mathrm{Spin}}
\newcommand\Pin{\mathrm{Pin}}
\newcommand\HSpin{\mathrm{HSpin}}
\newcommand\spin{\mathrm{Spin}}
\newcommand\SO{\mathrm{SO}}
\newcommand\Orth{\mathrm{O}}

\newcommand{\gm}{\GG_{\rmm}}
\newcommand{\ga}{\GG_{\rma}}

\newcommand{\dr}[1]{(\mspace{-3mu}(#1)\mspace{-3mu})}
\newcommand{\ds}[1]{[\mspace{-2mu}[#1]\mspace{-2mu}]}
\newcommand{\da}[1]{\pfister{\mspace{-3mu}{#1}\mspace{-3mu}}}

\renewcommand{\setminus}{\smallsetminus}

\begin{document}

\title{Essential dimension and algebraic stacks}

\author[Brosnan]{Patrick Brosnan$^\dagger$}

\author[Reichstein]{Zinovy Reichstein$^\dagger$}

\author[Vistoli]{Angelo Vistoli$^\ddagger$}

\address[Brosnan, Reichstein]{Department of Mathematics\\
The University of British Columbia\\
1984 Mathematics Road\\
Vancouver, B.C., Canada V6T 1Z2}

\address[Vistoli]{Scuola Normale Superiore\\Piazza dei Cavalieri 7\\
56126 Pisa\\Italy}

\email[Brosnan]{pbrosnan@math.ubc.ca}
\email[Reichstein]{reichst@math.ubc.ca}
\email[Vistoli]{angelo.vistoli@sns.it}


\begin{abstract} We define and study the essential dimension of an
  algebraic stack.  We compute the essential dimension of the stacks
  $\cM_{g,n}$ and $\overline{\cM}_{g,n}$ of smooth, or stable,
  $n$-pointed curves of genus~$g$.  We also prove a general lower
  bound for the essential dimension of algebraic groups with a
  non-trivial center. Using this, we find new exponential lower bounds
  for the essential dimension of spin groups and new formulas for
  the essential dimension of some finite $p$-groups.  Finally, we apply the
  lower bound for spin groups to the theory of the Witt ring of quadratic
  forms over a field $k$.
\end{abstract}
\thanks{$^\dagger$Supported in part by an NSERC discovery grant}
\thanks{$^\ddagger$Supported in part by the PRIN Project ``Geometria
sulle variet\`a algebriche'', financed by MIUR}

\maketitle

\tableofcontents

\section{Introduction}

Let $k$ be a field. We will write $\Fields_k$ for the category of
field extensions $K/k$.  Let $F\colon\Fields_k \arr\Sets$ be a covariant
functor.

\begin{definition} 
\label{def.merkurjev} Let $a\in F(L)$
  for $L$ an object of $\Fields_k$. A \emph{field of definition} for
  $a$ is an intermediate field $k \subseteq K \subseteq L$ such that
  $a$ is in the image of the induced function $F(K) \arr F(L)$.
  
  The \emph{essential dimension} $\ed a$ of $a$ (with respect to $L$)
  is the minimum of the transcendence degrees $\trdeg_{k}K$ taken over
  all fields of definition of $a$.

The \emph{essential dimension $\ed F$ of the functor} $F$ is the
supremum of $\ed a$ taken over all $a\in F(L)$ with $L$ in
$\Fields_k$.  \end{definition}

  Note that in Definition~\ref{def.merkurjev} the essential
  dimension of $a$ depends on the field $L$.   We write $\ed a$
  instead of $\ed(a,L)$ to simplify the notation. 

\begin{remark}
  If the functor $F$ is limit-preserving, a condition that is
  satisfied in all cases that interest us, every element $a \in F(L)$
  has a field of definition $K$ that is finitely generated over $k$,
  so $\ed a$ is finite. On the other hand, $\ed F$ may be infinite
  even in cases of interest (see for example Theorem~\ref{thm.curves}).
\end{remark}

\begin{example} \label{ex.edG}
Let $G$ be an algebraic group. Consider the Galois cohomology functor
$\H^1(*, G)$ sending $K$ to the set $\H^1(K, G)$ of 
isomorphism classes of $G$-torsors over $\Spec(K)$. The essential 
dimension of this functor is a numerical invariant of $G$, usually 
denoted by $\ed G$. Essential dimension was originally introduced 
(in~\cite{bur, reichstein}) and has since been extensively studied 
in this context; see, e.g.,~\cite{ry, kordonskii0, ledet, jly, bf1, lemire, 
cs, garibaldi}.
\end{example}

Definition~\ref{def.merkurjev} is due to A.~Merkurjev, 
as is the following observation; cf.~\cite[Proposition 1.17]{bf1}. 

\begin{example} \label{ex.ed-variety}
Let $X/k$ be a scheme of finite type over a field $k$, and let $F_X\colon\Fields_k  \arr\Sets$ 
denote the functor given by $K\arrto X(K)$.  Then $\ed F_X=\dim X$.
\end{example}

Note that the same is true if $X$ is an algebraic space (see 
Proposition~\ref{p.ed-variety}).

Many interesting naturally arising functors are not of the form
discussed in Examples~\ref{ex.edG} or~\ref{ex.ed-variety}.  One such
example is the functor $\operatorname{Curves}_{g,n}$ that sends $K$
into the set of isomorphism classes of $n$-pointed smooth algebraic
curves of genus $g$ over $K$. When $2g-2+n > 0$ this functor has a
well-known extension $\operatorname{\overline{Curves}}_{g,n}$ that
sends $K$ into the set of isomorphism classes of $n$-pointed stable
algebraic curves of genus $g$ over $K$.
Much of his paper was motivated by the following question. 

\begin{question} \label{q.curves} What are
  $\ed\operatorname{Curves}_{g,n}$ and
  $\ed\operatorname{\overline{Curves}}_{g,n}$?
\end{question}

Our starting point is the following definition.

\begin{definition} \label{def.ed-stack} Suppose $\cX$ is an algebraic
  stack over $k$.  The \emph{essential dimension} of $\cX$ is the
  essential dimension of the functor $F_{\cX}\colon\Fields_k \arr\Sets$ which
  sends a field $L/k$ to the isomorphism classes of objects in
  $\cX(L)$.  We write $\ed \cX$ for the \emph{essential dimension} of
  the stack $\cX$.
\end{definition}

Note that all of the examples above may be viewed as special cases
of~\ref{def.ed-stack}.  If $\cX$ is a scheme of finite type (or an
algebraic space), we recover Example~\ref{ex.ed-variety}. If $\cX =
\cB G$, the classifying stack of $G$ such that $\cB G(T)$ is the category of $G$-torsors on $T$, we recover Example~\ref{ex.edG}.  Finally, Question~\ref{q.curves} asks for the
values of $\ed\cM_{g,n}$ and $\ed\overline{\cM}_{g,n}$, where
$\cM_{g,n}$ and $\overline{\cM}_{g,n}$ are the stacks of $n$-pointed
smooth, or stable, algebraic curves of genus $g$ over a field $k$.
 
\begin{remark}
\label{r.abuse}
If $G$ is an algebraic group, we will often write $\ed G$ for $\ed\cB G$.
That is, we will write $\ed G$ for the essential dimension of the
stack $\cB G$ and not the essential dimension of the scheme
underlying $G$.  We do this to conform to the, now
standard, notation described in Example~\ref{ex.edG}. Of course, by
Example~\ref{ex.ed-variety}, the essential dimension of the underlying scheme is $\dim G$.
\end{remark}

In this paper we develop the theory of essential dimension 
for algebraic stacks. As a first application of this theory, we
give the following answer to Question~\ref{q.curves}.

\begin{theorem}
\label{thm.curves}Assume that the characteristic of $k$
  is $0$. Then
\[
   \ed\operatorname{Curves}_{g,n} = \ed\cM_{g,n}
   = 
   \begin{cases} 
   2         & \text{if }(g,n)=(0,0)\text{ or } (1,1);\\
   0         & \text{if }(g,n)=(0,1)\text{ or } (0,2);\\
   +\infty   & \text{if }(g,n)=(1,0);\\
   5         & \text{if }(g,n)=(2,0);\\
   3g-3 + n  & \text{otherwise}.
\end{cases}
\]

Moreover for $2g-2+n > 0$ we have $\ed\overline{\cM}_{g,n} = \ed \cM_{g,n}$.
\end{theorem}

Notice that $3g-3+n$ is the dimension of the moduli space $\rM_{g,n}$
in the stable range $2g-2+n > 0$ (or the dimension of the stack in all
cases); the dimension of the moduli space represents an obvious lower
bound for the essential dimension of a stack. The first
four cases are precisely the ones where a generic object in
$\cM_{g,n}$ has non-trivial automorphisms, and the case $(g,n) =
(1,0)$, is the only one where the automorphism group scheme 
of an object of $\cM_{g,n}$ is not affine.

Our stack-theoretic formalism turns out to be useful even 
for studying the essential dimension of algebraic groups  
in the classical setting of Example~\ref{ex.edG}.
Our key result in this direction is Theorem~\ref{thm2} below.

Let \begin{equation}
\label{e.extensions}
1\arr Z \arr G \arr Q\arr 1
\end{equation}
denote an extension of group schemes over a field $k$ with $Z$ central
and isomorphic to $\mu_n$ for some integer $n>1$.
For every extension $K$ of $k$ the 
sequence~\eqref{e.extensions} induces a connecting homomorphism
$\partial_K\colon \H^1(K,Q) \arr \H^2(K,Z)$.
We define $\ind(G,Z)$ as the maximal value of 
$\ind\bigl(\partial_K(t)\bigr)$ as
$K$ ranges over all field extensions of $k$ and $t$ ranges over all
torsors in $\H^1(K,Q)$. (Note that $\ind(G,Z)$ does not depend 
on the choice of the isomorphism $Z \simeq \mu_{n}$.)

\begin{theorem}
\label{thm2}
Let $G$ be an extension as in~\eqref{e.extensions}.  Assume that $n$
is a prime power. Then $\ed G \geq \ind(G,Z) - \dim Q$.
\end{theorem}

Let $G$ be a finite abstract group.  We write $\ed_k G$ for the
essential dimension of the constant group scheme $G_k$ over the field
$k$.  Let $\exp G$ denote the exponent of $G$ and let $\rC(G)$ denote
the center of $G$.
One of the main consequences of Theorem~\ref{thm2} is the following 
result about the essential dimension of finite $p$-groups. 

\begin{theorem}\label{t.p-groups}
Let $G$ be a $p$-group whose commutator $[G, G]$ is central and cyclic.
Then
   \[
   \ed_k G = \sqrt{|G/\rC(G)|} + \rank \, \rC(G) - 1 \, . 
   \]
for any base field $k$ of characteristic $\ne p$ which contains 
a primitive root of unity of degree $\exp(G)$.
\end{theorem}

Note that, with the above hypotheses, $|G/\rC(G)|$ is a complete
square.
In the case where $G$ is abelian we recover the identity $\ed(G) = \rank(G)$;
cf.~\cite[Example 7.4]{ry}. For most finite groups $G$ the best 
previously known lower bounds on $\ed(G)$ were of the form
\begin{equation} \label{e0.old-bound}
\ed(G) \ge \rank(A) \, ,
\end{equation}
where $A$ was taken to be an abelian subgroup $A$ of $G$ 
of maximal rank.  Theorem~\ref{t.p-groups} 
represents a substantial improvement over these bounds. 
For example, if $G$ is a non-abelian group of 
order $p^3$ and $k$ contains a primitive root of unity of 
degree $p^2$ then Theorem~\ref{t.p-groups} tells us that $\ed(G)= p$,
while~\eqref{e0.old-bound} yields only $\ed(G) \ge 2$. 

Theorem~\ref{t.p-groups} has a number of interesting consequences.
One of them is that $\ed(G) \ge p$ for any non-abelian $p$-group $G$;
see Corollary~\ref{cor2.p-groups}. Another is the following new bound on 
$\ed \Spin_n$.  Here by $\Spin_{n}$ we will mean the totally split
form of the the spin group in dimension~$n$ over a field $k$.

\begin{theorem} \label{thm.spin} Suppose $k$ is a field of 
characteristic $\neq 2$, and that $\sqrt{-1} \in k$. If 
$n$ is not divisible by $4$ then
\begin{align*}
   2^{\lfloor (n -1)/2 \rfloor} - \frac{n(n-1)}{2} &\le
   \ed \Spin_{n} \le  2^{\lfloor (n - 1)/2 \rfloor} \, . \\
\intertext{If $n$ is divisible by $4$ then}
   2^{\lfloor (n-1)/2 \rfloor} - \frac{n(n-1)}{2} + 1 &\le
   \ed \Spin_{n} \le  2^{\lfloor (n - 1)/2 \rfloor} + 1.
\end{align*}
\end{theorem}

The lower bound in this theorem was surprising to us because 
previously the best known lower bound was the following 
result due of V.~Chernousov and J.--P.~Serre~\cite{cs}.
\begin{equation} \label{e.old-spin}
\ed \Spin_n \geq \begin{cases} \lfloor n/2 \rfloor + 1 &\text{if $n
    \ge 7$
    and $n \equiv 1$, $0$ or $-1 \pmod{8}$} \\
  \lfloor n/2 \rfloor &\text{for all other $n \ge 11$.} \end{cases} 
\end{equation}
(The first line is due to B.~Youssin and the second author in the case
that $\chr k=0$~\cite{ry}.)  Moreover, in low dimensions,
M.~Rost~\cite{rost} (cf. also~\cite{garibaldi}) computed 
the following table of exact values:
\begin{center} \renewcommand{\arraystretch}{1.25}
\rule{0pt}{14pt}
\begin{tabular}{r@{${}={}$}lr@{${}={}$}lr@{${}={}$}lr@{${}={}$}l}
$\ed \Spin_3$ & $0$ &$ \ed \Spin_4$ & $0$ & $\ed \Spin_5$ & $0$ &
   $\ed \Spin_6$ & $0$ \\
$\ed \Spin_7$ & $4$ & $\ed \Spin_8$ & $5$ & $\ed \Spin_9$ & $5$ &
   $\ed \Spin_{10}$ & $4$\\
$\ed \Spin_{11}$ & $5$ & $\ed \Spin_{12}$ & $6$ & $\ed
\Spin_{13}$ & $6$ &
   $\ed \Spin_{14}$ & $7$.
\end{tabular}
\end{center}
Taken together these results seemed to suggest 
that $\ed\Spin_n$ should be a slowly increasing 
function of $n$ and gave no hint of its exponential growth.

Note that the computation of $\ed\Spin_n$ gives an example of a
split, simple, connected 
linear algebraic group whose essential dimension exceeds its
dimension.  (Note that for a simple adjoint 
group $G$, $\ed(G) \le \dim(G)$; cf. Example~\ref{ex.cgr1}.) 
It also gives an example of a split, semi-simple, connected
linear algebraic group $G$ with a central subgroup $Z$ such that 
$\ed G>\ed G/Z$.  This is because $\ed\SO_n = n - 1$ for $n \ge 3$;
cf.~\cite[Theorem 10.4]{reichstein}. 

Finally we follow a suggestion of A. Merkurjev and B. Totaro to
apply our results on $\ed\Spin_{n}$ to a problem in the theory 
of quadratic forms. Let $K$ be a field of characteristic 
different from $2$ containing a square root of $-1$, 
and let $\rW(K)$ be the Witt ring of $K$. Call $I(K)$ 
the augmentation ideal in $\rW(K)$; it is well known that 
if $q$ is a non-degenerate $n$-dimensional quadratic 
form whose class $[q]$ in $\rW(K)$ lies in $I^{a}(K)$, 
then $[q]$ can be expressed as the class a sum of $a$-fold Pfister forms. 
It is a natural to ask how many form are needed. 
When $a = 1$ or $a = 2$ is easy to see that $n$ Pfister 
forms suffice; see Proposition~\ref{p.PfisterEasy}. 
We prove the following result.

\begin{theorem} \label{t.pfister} 
Let $k$ be a field of characteristic different 
from $2$ and $n$ an even positive integer.
Then there is a field extension $K/k$ and
a class $[q] \in I^3(K)$ represented by an $n$-dimensional
quadratic form $q/K$ such that 
$[q]$ cannot be written as the 
sum of fewer than
   \[ \frac{2^{(n+4)/4} -n-2}{7}  \]
$3$-fold Pfister forms over $K$.
\end{theorem}

\subsection*{Description of contents}
The rest of this paper is structured as follows.

\S\ref{s.Gen} contains general results on essential dimension of
algebraic stacks, which are used systematically in the rest of the
paper.

\S\ref{s.quotient} contains a discussion of essential dimension of
quotient stacks; here we mostly rephrase known facts in our language.
At the end of the section, we show finiteness of the essential
dimension for a large class of algebraic stacks of finite type over a
field.
This class includes all Deligne--Mumford stacks and all 
quotient stacks of the form $[X/G]$ for $G$ a linear
algebraic group.

In \S\ref{s.generic} we prove Theorem~\ref{thm:generic} about
essential dimension of smooth integral \dm stacks satisfying an
appropriate separation hypothesis; it states that the essential
dimension of such a stack is the sum of its dimension and the
essential dimension of its generic gerbe. This somewhat surprising result
implies that the essential dimension of a non-empty open substack
equals the essential dimension of the stack. In particular, it proves
Theorem~\ref{thm.curves} in the cases where a general curve in
$\cM_{g,n}$ has no non-trivial automorphisms. It also brings 
into relief the important role played by gerbes in this theory.

Our main result on gerbes is Theorem~\ref{t.edGerbe}, stated in
\S\ref{s.gerbes} and proved in \S \ref{s.canonical} and
\S\ref{s.ed-gerbe}. It says that the essential dimension of a gerbe
banded by $\mu_{n}$, where $n$ is a prime power, equals the index of
its class in the Brauer group. Our proof is geometric: we link the
essential dimension of a gerbe banded by $\mu_{n}$ with the canonical
dimension of the associated Brauer--Severi variety, and use a result
of Karpenko on the canonical dimension of Brauer--Severi varieties of
prime-power index.

In \S \ref{s.ed-Mgn} we use Theorems~\ref{thm:generic}
and~\ref{t.edGerbe} to compute the essential dimensions of stacks of
hyperelliptic curves. This and some special arguments complete the
proof of Theorem~\ref{thm.curves}, except for the statement 
that $\ed \cM_{1,0} = +\infty$.

Theorem~\ref{t.edGerbe} is used again in \S\ref{s.extensions}, where
we prove Theorem~\ref{thm2}. 
The rest of the paper is dedicated to applications of this result.

In \S\ref{s.Tate} we complete the proof of Theorem~\ref{thm.curves} by
showing that $\ed \cM_{1,0} = +\infty$. This is achieved by applying
Theorem~\ref{thm2} to the group schemes of $l^{n}$-torsion points on
the Tate curves, where $l$ is a prime.

\S\ref{s.ed-pgroups} contains our results on $p$-groups.  
We also answer a question of Jensen, Ledet and Yui by 
giving an example of a finite group $G$ with a normal subgroup
$N$ such that $\ed(G/N) > \ed G$~\cite[p.204]{jly}.

Theorem~\ref{thm.spin} is proved in \S\ref{s.spinor}, along with
similar estimates for the essential dimensions of pin and 
half-spin groups.

Theorem~\ref{thm2} can also be applied to cyclic group over small fields. 
Little was known about the essential dimension of a cyclic group over
$\QQ$ until recently, when an important preprint~\cite{florence}
of M. Florence appeared, computing the essential dimension of 
a cyclic group of order $p^{m}$, where $p$ is a prime, 
over a field containing a primitive $p$-th root of $1$; 
this implies that $\ed_{\QQ} (\ZZ/p^m)\geq p^{m-1}$.  
In \S\ref{s.cyclic} we we recover this result as a consequence of
Theorem~\ref{thm2} by making use of the Brauer-Rowen algebra, 
an idea we learned from~\cite{florence}. 
As a corollary of Florence's theorem
we prove a particular case of a conjecture of Ledet~\cite[Section 3]{ledet}, 
relating the essential dimensions of the cyclic group $C_n$ and
the dihedral group $D_{n}$ ($n$ odd). We show that $C_n$ 
and $D_n$ have the same essential dimension if $n$ is 
a prime power and $k$ contains a primitive $p$-th root of $1$.

Section \S\ref{s.arason} contains our application 
of the results on $\ed\Spin_{n}$ to the theory 
of quadratic forms.  In particular, we prove 
Theorem~\ref{t.pfister}.

\begin{remark}
  One interesting issue that we do not address in this paper is the
  subject of the essential dimension at a prime $p$.  For groups, this
  is defined in~\cite[Definition 6.3]{ry}.  There is an obvious
  generalization to stacks which we leave to the reader to formulate.
  We hope that the main results of this paper would remain valid for
  the concept of essential dimension at a prime $p$ (with certain
  obvious alterations).  In particular, we think it is very likely
  that (in the notation of ~\cite{ry}), $\ed(G;p)$ is given by the
  formula in Theorem~\ref{t.p-groups} and $\ed(\Spin_n;2)$ is bounded
  by the formulas in Theorem~\ref{thm.spin}.  However, we have not
  checked this in detail.
\end{remark}

\subsection*{Notation}
  In the paper, a \emph{variety} over a field $k$ will be a
  geometrically integral separated scheme of finite type over $k$.
  Cohomology groups $\H^i(T,\cF)$ will be taken with respect to the
  fppf topology unless otherwise specified.
  
  As explained in Remark~\ref{r.abuse}, we will write $\ed G$ for
  $\ed\cB G$ and use these notations interchangeably.  The reader may
  notice that we prefer to write $\ed\cB G$ earlier in the paper where
  we are working in a general stack-theoretic setting and $\ed G$
  towards the end where we are primarily concerned with essential
  dimensions of algebraic groups.
   
  We write $\mu_n$ for the groups scheme of $n$-th roots of unity.  If
  $k$ is a field, we write $\zeta_n$ for a primitive $n$-th root of unity in
  the algebraic closure of $k$.  (Using the axiom of choice, we choose
  one once and for all.)
  For typographical reasons, we sometimes write $C_n$ for the cyclic
  group $\ZZ/n$.

\subsection*{Acknowledgments}
We would like to thank the Banff International Research Station in
Banff, Alberta (BIRS) for providing the inspiring meeting place where
this work was started.  We are grateful to K.~Behrend, C.-L. Chai,
D.~Edidin, N.~Fakhruddin, A.~Merkurjev, B.~Noohi, G.~Pappas,
D.~Saltman and B.~Totaro for very useful conversations.  We would also
like to thank M.~Florence for sending us a copy of his preprint on the
essential dimension of $\ZZ/p^n$ over $\QQ(\mu_p)$.

\section{Generalities}
\label{s.Gen}

We begin by reformulating Definition~\ref{def.merkurjev} in the
language of fibered categories.  For this notion, we refer the reader
to the Definition 3.1 of~\cite{v1}.

For a field $k$, let $\zp_k\eqdef \Fields_k{\op}$. We study
categories $\cX$ which are fibered over $\zp_k$; these are a
generalization of functors from $\zp_k$ to $\sets$.  (Clearly $\cX$ is
fibered over $\zp_k$ if and only if $\cX$ is \emph{cofibered over
  $\Fields_k$} (\cite[6.10]{SGA1b}), but we prefer to work with fibered
categories.)

\begin{definition}
If $\xi$ is an object of $\cX(K)$, where $K$ is an extension of $k$, 
a \emph{field of definition} of $\xi$ is an intermediate field 
$k \subseteq F \subseteq K$, such that $\xi$ is in the essential 
image of the pullback functor $\cX(F) \arr \cX(K)$.
\end{definition}

\begin{definition}
\label{d.EdFibered}
Let $\cX$ be a category fibered over $\points_k$. 
If $K$ is an extension of $k$ and $\xi$ is an object 
of $\cX(K)$, the \emph{essential dimension of $\xi$}, written $\ed\xi$, is 
the least transcendence degree over $k$ of a field of definition of $\xi$.

The \emph{essential dimension of $\cX$}, denoted by $\ed\cX$, 
is the supremum of the essential dimension of all objects 
$\xi$ in $\cX(K)$ for all extensions $K$ of $k$.
\end{definition}  

These notions are obviously relative to the base field $k$.  (See
Remark~\ref{r.BaseField}.)  We will write $\ed(\xi/k)$ (resp.
$\ed(\cX/k)$) when we need to be specific about the dependence on the
base field.  

Note that $\ed\cX$ takes values in the range
$\{\pm\infty\}\cup\ZZ_{\geq 0}$, with $-\infty$ occurring if and only if
  $\cX$ is empty.

\begin{remark}
\label{r.functors}
  With every functor $F\colon\Fields_k \arr\sets$, one can canonically
  associate a category $\cX_F$ fibered over $k$ (see~\cite[Proposition
  3.26]{v1}).  It is an easy exercise in unravelling the
  definitions to see that $\ed\cX_F$ as defined in~\ref{d.EdFibered}
  is equal to $\ed F$ as defined in~\ref{def.merkurjev}.
  
  Furthermore, given a fibered category $\cX \arr \Aff_{k}$ we get a
  functor
   \[
   \overline{\cX}\colon \Fields_{k} \arr \Sets
   \]
sending a field $K$ into the set of isomorphism classes in $\cX(\spec K)$. It also straightforward to see that $\ed\cX$ equals $\ed\overline{\cX}$ as  defined in~\ref{def.merkurjev}.
\end{remark}

\begin{remark} 
\label{r.BaseField}
Let $L$ be an extension of a field $K$ of
  transcendence degree~$d$, and let $\cX \arr \zp_{L}$ be a fibered
  category. Then the composite of $\cX \arr \zp_{L}$ with the obvious
  functor $\zp_{L} \arr \zp_{K}$ makes $\cX$ into a category fibered
  over $\zp_{K}$. We have
   \[
   \ed(\cX/K) = \ed(\cX/L) + d.
   \]
The idea is that if $F$ is an extension of $K$ and $\xi$ is an object of
$\cX(F)$, the image of $\xi$ in $\zp_{L}$ defines an embedding of $L$
in $F$; and every field of definition of $\xi$ over $K$ must contain
$L$.
\end{remark}

\begin{para} If $\xi$ is an object in $\cX(L)$ and $k\subset F\subset
  L$, then we say that $\xi$ \emph{descends to $F$} if $F$ is a field
  of definition for $\xi$.   The map $\Spec L\to \Spec F$ is then
  called a \emph{compression} of $\xi$.   Note that the identity
  morphism $\id_{\Spec L}$ is always a compression.
  
  A compression $\Spec L \arr \Spec F$ is called a \emph{deflation} 
  if  $\trdeg_k F<\trdeg_k L$. If a deflation exists, we will
  say that $\xi$ is \emph{deflatable}. If $\xi$ is not
  deflatable, then it is called \emph{undeflatable}. To show that
  $\ed\cX\geq d$, it suffices to exhibit an undeflatable object
  $\xi\in\cX(L)$ with $\trdeg_k L=d$.

  We will call an undeflatable object $\xi\in\cX(L)$ \emph{maximally
    undeflatable} if $\ed\xi=\ed\cX$. Clearly, for $d\in [0,\infty)$,
  we have $\ed (\cX/k) =d$ if and only if there is a maximal undeflatable
  object $\xi\in\cX(L)$ with $\trdeg_k L=d$.
\end{para}

\begin{para}
  \label{p.BaseCategory} For a scheme $S$, we follow Laumon \&
  Moret-Bailly~\cite{LMB} in letting $\Aff_S$ denote the category of
  affine schemes over $S$.  If $S=\Spec k$, then this category, which
  we denote by $\Aff_k$, is equivalent to the category opposite to the category
  $k$-$\Alg$ of $k$-algebras.  We equip $\Aff_S$ with the \'etale
  topology, and, by default, all notions of sheaves and stacks
  involving $\Aff_S$ are with respect to this topology.
\end{para}

\begin{para}
\label{p.stacks}
A \emph{stack} over a scheme $S$ will mean a stack over $\Aff_S$.  That
is, a stack over $S$ is a category $\cX$ fibered over $\Aff_S$
satisfying Definition 4.6 of~\cite{v1}.  If $\cX$ is a
category fibered over $\Aff_k$, then the restriction $\widetilde{\cX}$ of
$\cX$ to $\zp_k$ via the obvious functor $\zp_k \arr\Aff_k$ is a
category fibered over $\zp_k$.  We write $\ed\cX\eqdef \ed\widetilde\cX$.  
This defines the notion of the essential dimension of a stack.
\end{para}

\begin{para}
  \label{p.algstacks}
  We use~\cite[Definition 4.1]{LMB} as our definition of an algebraic
  stack.  That is, by an \emph{algebraic stack} over a scheme $S$, we
  will mean a stack $\cX$ in groupoids over $\Aff_S$ satisfying
\begin{enumerate}
\item The diagonal morphism $\Delta\colon\cX \arr\cX\times_S\cX$ is
  representable, separated and quasi-compact,
\item there is an algebraic space $X$ over $S$ and a morphism
$X \arr\cX$ which is smooth and surjective.
\end{enumerate}
\end{para}

We will make heavy use of the notion of gerbe. Let us recall that
a category $\cX$ fibered in groupoids over the category $\Aff_K$ is
an \emph{fppf gerbe} if the following conditions are satisfied.

\begin{enumerate}

\item $\cX$ is a stack with respect to the fppf topology.

\item There exists a field extension $K'$ of $K$ such that $\cX(\spec
  K')$ is not empty.

\item Given an affine scheme $S$ over $K$ and two objects $\xi$ and
  $\eta$ in $\cX(S)$, there exists an fppf cover $\{S_{i} \arr S\}$
  such that the pullbacks $\xi_{S_{i}}$ and $\eta_{S_{i}}$ are
  isomorphic in $\cX(S_{i})$ for all $i$.

\end{enumerate}

A gerbe is called \emph{neutral} if $\cX(\spec K)$ is not empty.

We have the following easy observation.

\begin{proposition}
\label{p.reduced}
Let $\cX$ be an algebraic stack over a field $k$, and let $\cX_{\red}$
denote the reduced substack~\cite[Lemma 4.10]{LMB}.  Then $\ed\cX_{\red}=\ed\cX$.
\end{proposition}
\begin{proof}
  For every field $K$ over $k$, the morphism $\cX_{\red}\arr\cX$
  induces an equivalence of categories $\cX_{\red}(K)\to \cX(K)$.
\end{proof}

A category $\cX$ fibered over $\points_k$ (resp. over $\Aff_k$) is
\emph{limit preserving} if, whenever $K=\colim K_i$ is a filtered
direct limit of fields (resp. $k$-algebras), $\colim X(\Spec K_i) \arr
X(K)$ is an equivalence of categories (see~\cite[p. 167]{artin}).
Note that an algebraic stack (viewed as a category fibered over
$\Aff_k$) is limit preserving~\cite[Proposition 4.18]{LMB}.  The
property of being limit preserving provides the most basic instance of
finiteness of essential dimension.

\begin{proposition}
If $\cX$ is a limit-preserving category fibered over $\zp_k$, then 
any object $\xi$ of $\cX$ has finite essential dimension.
\end{proposition}
\begin{proof}
Let $K$ be a field.  We can write it as a filtered direct limit
$K=\colim_I K_i$ of all of its subfields $K_i$ of finite transcendence degree.
Since $\cX$ is limit preserving, every object $\xi$ of $X(K)$ must be
in the essential image of $\cX(K_i) \arr \cX(K)$ for some $i\in I$.
Therefore $\ed\xi<\infty$.  
\end{proof}

Let $K$ be an extension of $k$; there is a tautological functor
$\zp_{K} \arr \zp_{k}$. We denote by $\cX_{K}$ the pullback of $\cX$
to $K$. This means the following. An object of $\cX_{K}$ is a pair
$(\xi,E)$, where $E$ is an extension of $K$, and $\xi$ is an object of
$\cX(E)$ mapping to $E$ in $\points_{k}$. An arrow from $(\xi,E)$ to
$(\xi',E')$ is an embedding $E' \into E$ (corresponding to an arrow $E
\arr E'$ in $\points_{k}$) preserving $K$, and an arrow $\xi' \arr \xi$
in $\cX$ mapping to this embedding. The category $\cX_{K}$ is fibered
over $\points_{K}$.

\begin{proposition} \label{p.extend-field}
Suppose that $k$ is a field, $K$ is an extension of $k$ and $L$ is an extension of $K$.
\begin{enumeratea}
\item Let $\cX$ be a category fibered over $\zp_{k}$. If $\xi$ is 
an object of $\cX(K)$ and $\xi_{L}$ denotes its pullback in $\cX(L)$, then
   \[
   \ed \xi_{L} \leq \ed \xi.
   \]

\item Let $\cX$ be a limit-preserving fibered category over $\Aff_{k}$. If $\xi$ is an object of $\cX$ and $L$ is contained in a purely transcendental extension of $K$, then
   \[
   \ed \xi_{L} = \ed\xi.
   \]

\end{enumeratea}
\end{proposition}

\begin{proof}
Part~(a) is obvious.

For part~(b), let $E$ be a purely transcendental extension of $k$ containing $K$: then $\ed\xi_{E} \leq \ed\xi_{L} \leq \ed\xi$ by part~(a), so it is enough to prove that $\ed\xi_{E} = \ed\xi$. So we can assume that $L$ is purely transcendental over $k$.

We need to prove the inequality $\ed\xi_{L} \geq \ed\xi$. Let $B$ be a
transcendence basis of $L$ over $K$. Let $E \subseteq L$ be a finitely
generated subfield with $\trdeg_{k}E = \ed\xi_{L}$, with an object
$\xi_{E}$ whose image in $\cX(L)$ is isomorphic to $\xi_{L}$; then if
$E \subseteq L'\subseteq L$ is an intermediate field and $\xi_{L'}$ is
the image of $\xi_{E}$ in $\cX(L')$, we have $\ed\xi_{L'} =
\ed\xi_{L}$. Since $L$ is the directed limit of subfields of the form
$K'(S)$, where $K' \subseteq K$ is a field of definition of $\xi$
which is finitely generated over $k$ and $S \subseteq B$ is a finite
subset, after enlarging $L'$ we can find such and $L'$ of the form
$K'(S)$. Let $\xi_{K'}$ be an object of $\cX(K')$ whose image in
$\cX(K)$ is isomorphic to $\xi$. The image of $\xi_{K'}$ in $\xi_{L'}$
is not necessarily isomorphic to $\xi_{L'}$, but it will become so
after enlarging $K'$ and $S$. Since we have $\ed\xi_{K'} \geq \ed\xi$
and $\ed \xi_{L'} = \ed\xi_{L}$, we may substitute $K'$ and $L'$ for
$K$ and $L$ and assume that $K$ and $L$ are finitely generated over
$k$.

We may also assume that $K$ is infinite, because otherwise $K$ is
finite over $k$, we have $\ed\xi = 0$ and the inequality is
obvious. Again because $\cX$ is limit-preserving, there will be an
affine integral scheme $U$ of finite type over $k$, with quotient
field $K$, and an object $\xi_{U}$ in $\cX(U)$, whose image in
$\cX(K)$ is isomorphic to $\xi$. Let $\{x_{1}, \dots, x_{n}\}$ be a
transcendence basis for $L$ over $K$. There will exist an open affine
subscheme $V$ of $U \times \AA_{k}^{n}$ and a dominant morphism $V
\arr W$ onto an integral affine scheme $W$ which is of finite type of
dimension~$\ed\xi_{L}$ over $k$, together with an object $\xi_{W}$
whose pullback in $\cX(V)$ is isomorphic to the pullback of $\xi_{U}$
along the first projection $\pr_{1}\colon W \arr V$. Since the
fraction field of $U$ is infinite, there will exist a non-empty open
subscheme $U' \subseteq U$ and a section $U' \arr W$ of $\pr_{1}\colon
W \arr V$. From this we we see that the restriction $\xi_{U'}$ of
$\xi_{U}$ to $U'$ is isomorphic to the pullback of $\xi_{W}$ to
$\cX(U')$. If $V'$ denotes the closure of the image of $U'$ into $V$,
we get an object $\xi_{V'}$ of $\cX(V')$ whose image in $\cX(U')$ is
isomorphic to $\xi_{U'}$. Hence $k(V') \subseteq k(U') = K$ is a field
of definition of $\xi$; since $\dim V' \leq \dim V = \ed\xi_{L}$ we
conclude that $\ed\xi \leq \ed\xi_{L}$.
\end{proof}

The following observation is a variant of~\cite[Proposition 1.5]{bf1}.
We will use it repeatedly in the sequel.

\begin{proposition}
\label{p.extensions}
Let $\cX$ be a category fibered over $\zp_k$, and let $K$ 
be an extension of
$k$. Then $\ed(\cX_{K}/K) \leq \ed(\cX/k)$.
\end{proposition}
\begin{proof}
If $L/K$ is a field extension,
then the natural morphism $\cX_K(L) \arr \cX(L)$ is an equivalence.
Suppose than $M/k$ is a field of definition for an object $\xi$ in
$\cX (L)$. Then any field $N$ containing both $M$ and $K$ is a field of
definition for $\xi$.  Since there is a field $N/K$ with $\trdeg_K
N\leq \trdeg_k M$ containing both $M$ and $K$ (any composite of $M$
and $K$), we can find a field of definition for $\xi$ as an object in
$\cX_K$ of transcendence degree $\leq \ed(\cX/k)$. 
\end{proof}

\begin{remark}\label{r.extensions}
  The proof shows the following: if $L$ is an extension of $K$ and
  $\xi$ is an object in some $\cX_{K}(L)$, call $\eta$ the image of
  $\xi$ in $\cX$. Then $\ed \xi \leq \ed \eta$.
\end{remark}

In some cases of interest we can arrange for equality
Proposition~\ref{p.extensions}.

\begin{proposition}
\label{p.extensions2}
  Let $\cX$ be a limit-preserving category fibered over $\Aff_k$.  Suppose
one of the following conditions holds:
\begin{enumerate}
\item $k$ is algebraically closed.
\item $K/k$ is purely transcendental and $k$ is infinite,
\end{enumerate}
Then $\ed(\cX_K/K)=\ed(\cX/k)$.
\end{proposition}
\begin{proof}
It is easy to see that, since $\cX$ is limit-preserving, we can 
assume that $K$ is finitely generated over $k$.   Thus, $\trdeg_{k} K<\infty$.  

So pick $\xi\in\cX(l)$ an undeflatable object for some  field extension $l/k$ with
$\trdeg_k l=n<\infty$.  Again, since $\cX$ is limit-preserving, we
can assume that $l$ is finitely generated over $k$.
Set $L\eqdef l\otimes_k K$.  Note that $L$ is a field of transcendence
degree $n$ 
under either hypotheses (1) or (2).    Write $\eta$ for the
restriction of $\xi$ to $L$ via the obvious map $\Spec L \arr \Spec l$.
We claim that, in either case (1) or (2), $\eta$ is undeflatable
over $K$.

To show this, assume that $\eta$ is deflatable over $K$.  Then there
is a  intermediate field $K\subset R\subsetneq L$
and an object
$\gamma\in\cX(R)$ such that the restriction of $\gamma$ to $L$ is
$\eta$.  Moreover $\trdeg_K R<\trdeg_K L$.
Pick affine schemes $U$ and $V$ of finite type over $k$ such that
$k(U)=l$ and $k(V)=K$. In case (2) we can and will assume that
$V=\AA^m$.  Since $\cX$ is limit-preserving, we can find
an affine scheme $Z$ of finite type over $k$ such that $k(Z)=R$ and
$\gamma$ is the restriction of some object $\widetilde{\gamma}\in\cX(Z)$.  
Shrinking $Z$ if necessary, we can assume that the rational
map from $Z$ to $V$ inducing the inclusion of $K$ into $R$ is a morphism.
We can also find a Zariski dense open $W\subset U\times V$ and a $V$-morphism 
$f\colon W \arr Z$ inducing the inclusion of $R$ into $L$.  

For each point $v\in V$,
write $W_v$ (resp. $Z_v$) for the fiber of the map $Z \arr V$ (resp. the
map $p_2\colon W \arr V$).  Since $R\subsetneq L$, $\dim W>\dim
Z$. It follows that there is a Zariski dense open subscheme $M$ of $V$ 
such that $\dim W_v > \dim Z_v$ for all $v\in M$. 

Suppose $v$ is a closed point in $M$ such that $k(v)=k$. Then $\xi$ is the
restriction of $\gamma$ to $l$ via the map $\Spec l \arr W_v \arr
Z_v$.  It follows that $\xi$ is deflatable over $k$, which is a
contradiction.

Note that, in either case (1) or case (2), $M(k)\neq\emptyset$.
Therefore the contradiction is always obtained.  It follows that
$\eta$ is undeflatable over $L$.

Now, we are free to pick $\xi\in\cX(l)$ with $\ed (\xi/k) =\ed (\cX/k)$.
Then, in either case (1) or (2), we have that $\ed (\cX/k)=\trdeg_k l
=\trdeg_K L=\ed (\cX/L)$.  This completes the proof of the statement.
\end{proof}

We will need the following generalization of Example~\ref{ex.ed-variety}.

\begin{proposition}
\label{p.ed-variety}
The essential dimension of an algebraic space locally 
of finite type over $k$ equals its dimension.
\end{proposition}
\begin{proof}
 Indeed, in this case  $X$ has a stratification by schemes $X_i$.  
Any $K$-point $\eta\colon\Spec K \arr X$ must land in one of the $X_i$.  
Thus $\ed X=\max\ed X_i=\dim X$.
\end{proof}

\begin{para}
  If $X$ is an algebraic space over an algebraic space $S$, then the
  category of arrows $T \arr X$ where $T$ is an object in $\Aff_S$ is
  fibered over $\Aff_S$.  It is equivalent to the fibered category
  $\cX_{\rmh_X}$ arising from the functor $\rmh_X\colon\Aff_S \arr\Sets$ given by
  $\rmh_X(T)=\Mor_S (T,X)$ via~\cite[Proposition 3.26]{v1}.

  A category $\cX$ fibered over $\Aff_S$ is said to be
  \emph{representable} by an algebraic space if there is an algebraic
  space $X$ over $S$ and an equivalence of categories between $\cX$
  and $\cX_{\rmh_X}$.  We will follow the standard practice of
  identifying an algebraic space $X$ with its corresponding
  representable stack $\cX_{\rmh_X}$.  This is permissible by Yoneda's
  lemma.

A morphism $f\colon\cX \arr\cY$ of categories fibered over $S$ is said to be
representable if, for every algebraic space $T$ over $\cY$, the fiber
product $\cX\times_{\cY} T$ is representable as a category fibered
over $\Aff_T$.

Let $d$ be an integer, and let $k$ be a field.  A morphism $f\colon\cX
\arr\cY$ of categories fibered over $\Aff_k$ is said to be
\emph{representable of fiber dimension at most $d$} if, for every map
$T \arr \cY$ from an algebraic space, the fibered product
$\cX\times_{\cY} T$ is an algebraic space locally of finite type over
$T$ with fibers of relative dimension $\leq d$.
\end{para}

\begin{proposition}
\label{p.RelDim}
  Let $\cX$ and $\cY$ be fibered categories over $k$. Let $d$ be an
  integer, and assume that there exists a morphism $\cX \arr
  \cY$ that is represented by morphisms locally of finite algebraic
  spaces, with fiber dimension at most $d$. Then
$\ed(\cX/k) \leq \ed(\cY/k) + d$.
\end{proposition}

\begin{proof}
  Let $K$ be a field over $k$ and let $x\colon\Spec K \arr \cX$ be an
  object of $\cX(K)$.  Then $f\circ x\colon\Spec K \arr \cY$ is an
  object of $\cY(K)$, and we can find a field $L$ with a morphism
  $y\colon\Spec L \arr Y$ such that $k\subset L \subset K$,
  $\trdeg_k L\leq \ed \cY$ and the following diagram commutes.
\begin{equation*}
\xymatrix{
\Spec K   \ar[r]^-{x}\ar[d]  & \cX\ar[d]^{f}\\
\Spec L   \ar[r]^-{y}        & \cY      \\
}
\end{equation*}

Let $\cX_L\eqdef \cX\times_{\cY}\Spec L$.  By the hypothesis, $\cX_L$ is an
algebraic space, locally of finite type over $L$ and of relative
dimension at most $d$.  By the commutativity of the above diagram, the
morphism $x\colon\Spec K \arr \cX$ factors through $\cX_L$.  Let $p$ denote
the image of $x$ is $\cX_L$.  Since $\cX_L$ has dimension at most
$d$, we have $\trdeg_k k(p)\leq d$.  Therefore $\trdeg_k k(p)\leq
\ed\cY + d$.  Since $x$ factors through $\Spec k(p)$ the result
follows.
\end{proof}

\begin{remark}
\label{r.RelDim}
  The proof of Proposition~\ref{p.RelDim} clearly shows the following: 
For any field $K/k$ and any $\xi\in\cX(K)$, $\ed \xi\leq \ed f(\xi)+d$.
\end{remark}

The following simple observation will be used often in this paper.  

\begin{proposition}
\label{p.observation}
  Let $U$ be an integral algebraic space locally of finite type over
  $k$ with function field $K\eqdef k(U)$, and let 
$f\colon\cX \arr U$ 
be a stack over $U$.   Let $\cX_K$ denote the pullback of $\cX$ to
$\Spec K$.  Then 
\[
\ed \cX \geq \ed (\cX_K/K) + \dim U.
\]
\end{proposition}
\begin{proof}
  If $\Spec L \arr \cX_K$ is maximally undeflatable over $K$, then
  the morphism $\Spec L \arr \cX$ obtained by composing with the
  canonical morphism $\cX_K \arr \cX$ is maximally undeflatable
  over $k$.
\end{proof}

Let $\cX$ be a locally noetherian stack over a field $k$ with
presentation $P\colon X \arr\cX$.  Recall that the \emph{dimension
of $\cX$ at a point $\xi\colon\Spec K \arr \cX$} is given by
$\dim_x(X)-\dim_x P$ where $x$ is an arbitrary point of $X$ lying over
$\xi$~\cite[(11.14)]{LMB}. Let $\cY$ be stack-theoretic closure of the
image of $\xi$; that is, the intersection of all the closed substacks
$\cY_{i}$ such that $\xi^{-1}(\cY_{i}) = \spec K$. The morphism $\xi$
factors uniquely through $\cY \subseteq \cX$. We defined \emph{the
  dimension of the point $\xi$} to be the dimension of the stack $\cY$
at the point $\spec K \arr \cY$.

\begin{proposition}
\label{p.RelDim2}
  Let $\cX \arr \cY$ be a morphism of algebraic stacks over a field $k$.  Let 
$K/k$ be a field extension and let $y\colon\Spec K \arr \cY$ be a
point of 
dimension $d\in\ZZ$.  Let 
$\cX_K\eqdef \cX\times_{\cY}\Spec K$.  Then
\begin{equation*}
  \ed (\cX_K/K)\leq \ed (\cX/k) - d
\end{equation*}
\end{proposition}

\begin{proof}
By~\cite[Theorem 11.5]{LMB}, $\cY$ is the disjoint union of a finite
family of locally closed, reduced substacks $\cY_i$ such that each
$\cY_i$ is an fppf gerbe over an algebraic space $X_i$ with structural
morphism $A_i\colon\cY_i \arr Y_i$.  We can therefore replace $\cY$ by one of
the $\cY_i$ and assume that $\cY$ is an fppf gerbe over an algebraic
space $Y$.  Without loss of generality, we can assume that $Y$ is an
integral affine scheme of finite type over $k$.  

Let $p$ be the image of $\xi$ in $Y$.  
Since $\cY$ is limit-preserving, we can find an integral affine scheme
$U$ equipped with a morphism $i\colon U \arr \cY$ and a dominant morphism
$j\colon\Spec K \arr U$ such that $y$ is equivalent to $i\circ j$. 
We can also assume that the composition $U \arr
\cY \arr Y$ is dominant. 

Since $\cY$ is a gerbe over $Y$, it follows that $U \arr \cY$ is
representable of fiber dimension at most $\dim U-d$.  Now, form the
following diagram with Cartesian squares.
$$
\xymatrix{
  \cX_K\ar[r]\ar[d] & \Spec K\ar[d]\\
  \cX_U\ar[r]\ar[d] & \ar[d]U\\
  \cX\ar[r]         & \cY
}
$$
Since the vertical maps in the lower square are representable of fiber
dimension at most $\dim U-d$, 
\begin{align*}
\ed (X_K/K) &\leq \ed (\cX_{k(U)}/k(U))\\
            &\leq \ed\cX_U -\dim U\\
            &\leq \ed \cX  +\dim U -d + \dim U\\
            &\leq \ed\cX -d.\qedhere
\end{align*}
\end{proof}

In general, the inequality of Proposition~\ref{p.RelDim} only goes in
one direction. However, in important special cases we can obtain an
inequality in the reverse direction.

\begin{para} 
  \label{d.isotropic} We will say that a morphism $f\colon\cX \arr\cY$
  of categories fibered over $\points_k$ is \emph{isotropic} if for
  every extension $K$ of $k$ and every object $\eta$ of $\cY(K)$ there
  exists an object $\xi$ of $\cX(K)$ such that $f(\xi)$ is isomorphic
  to $\eta$.
\end{para}

\begin{proposition}
\label{p.isotropic}
  Let $f\colon\cX \arr\cY$ be an isotropic morphism of categories fibered
  over $\points_k$.  Then $\ed\cX\geq \ed\cY$.
\end{proposition}
\begin{proof}
  Let $K$ be an extension of $k$ and $\eta$ an object of $\cY(K)$. If
  $\xi$ is an object of $\cX(K)$ such that $f(\xi)$ is isomorphic to
  $\eta$, then a field of definition for $\xi$ is also a field of
  definition for $\eta$.
\end{proof}

\begin{remark}
  One obvious example of an isotropic morphism is the total space of a
  vector bundle over a Deligne-Mumford stack.  Any open substack of a
  vector bundle which is dense in every fiber is also isotropic.
\end{remark}

We will use the following proposition.

\begin{proposition}
\label{r.products}
  Let $\cX$ and $\cY$ be categories fibered over  $\points_k$.  Then 
$$
\ed (\cX\times_{\points_k} \cY) \leq \ed\cX + \ed\cY.
$$
\end{proposition}
\begin{proof}
  This is equivalent to Lemma 1.11 of~\cite{bf1}. The proof is
  immediate: if $(\xi, \eta)$ is an object in some $(\cX \times
  \cY)(K)$, then $k \subseteq F \subseteq K$ is a field of definition
  for $\xi$ with $\trdeg_{k} F \leq \ed \cX$ and $k \subseteq L
  \subseteq K$ is a field of definition for $\eta$ with $\trdeg_{k} L
  \leq \ed \cY$, then the subfield of $K$ generated by $F$ and $L$ is
  a field of definition for $(\xi,\eta)$, of transcendence degree at
  most $\trdeg_{k} F + \trdeg_{k} L \leq \ed\cX + \ed\cY$.
\end{proof}

\begin{remark}
The inequality in Proposition~\ref{r.products} is often strict.  For
example, let $k=\CC$, $\cX=\cB\mu_{2}$ and $\cY=\cB\mu_{3}$.  Then
$\cX\times_{\points_k} \cY=\cB \mu_{6}$.  However, we have $\ed\mu_{n} = 1$
for all integers $n>1$ by~\cite[Theorem 5.3]{bur}.
\end{remark}

\section{Quotient stacks and Finiteness}
\label{s.quotient}

Suppose a linear algebraic group $G$ is acting on an algebraic space
$X$ over a field $k$. We shall write $[X/G]$ for the quotient stack
$[X/G]$.  The functor $F_{[X/G]}$ associates to a field $K/k$ the set
isomorphism classes of diagrams
\begin{equation} \label{e.functor}
 \xymatrix{
T \ar@{->}[r]^{\psi} \ar@{->}[d]^{\pi} &  X \cr 
\Spec(K) &  } 
\end{equation} 
where $\pi$ is a $G$-torsor and $\psi$ is a $G$-equivariant map.

If $G$ is an algebraic group over $k$, then $\cB G\eqdef [\Spec k/G]$.
The functor $F_{\cB G}$ is equal to the functor $K\arrto \H^1 (K,G)$
sending $K$ to the isomorphism classes of $G$-torsors over $K$.

\begin{remark}
  As noted in the introduction, for $G$ an algebraic group, the
  essential dimension $\ed \cB G$ is equal to the essential dimension
  of Example~\ref{ex.edG} classically denoted by $\ed G$.  To prevent
  confusion, we remind the reader that we will use the notations
  $\ed\cB G$ and $\ed G$ interchangeably (as in Remark~\ref{r.abuse}).
\end{remark}

\begin{proposition}\label{p.RelDim-quotients}
Let $G \arr \spec K$ be an algebraic group acting on an algebraic space $X$ over $K$ and let $H$ be a closed subgroup of $G$. Then
    \[
   \ed{[X/H]} \leq \ed{[X/G]} + \dim G - \dim H.
   \]
\end{proposition}

\begin{proof}
The obvious morphism $[X/H] \arr [X/G]$ has fibers of dimension $\dim G - \dim H$, so this is a consequence of Proposition~\ref{p.RelDim}.
\end{proof}

\begin{lemma} \label{lem2.1} Suppose a linear algebraic group $H$ is
  acting on an algebraic space $X$. If $H$ is a subgroup of another
  linear algebraic group $G$ then the quotient stacks $[X/H]$ and
  $[X*_H G/G]$ are isomorphic.
\end{lemma}

\begin{proof} Here $X*_H G$ is the quotient of $X\times G$ by the $H$ action 
given by $h(x,g)=(xh^{-1},hg)$. 
This fact is standard but it is as easy to
prove it as it is to look for a reference. 

Note that, when $H$ acts
freely on $X$, the quotients $X/H$ and $X*_H G/G$ are both algebraic spaces

Let $E$ be an object in $[X/H]$, i.e, an $H$-torsor over a $k$-scheme
$S$ equipped with an $H$-equivariant map to $X$.  We associate to $E$
the $G$-torsor $E*_H G$ equipped with its natural morphism to the
algebraic space $X*_H G$.

On the other hand, suppose $F$ is an object in $[X*_H G/G]$, i.e., a
$G$-torsor over a $k$-scheme $S$ equipped with a $G$-equivariant map
to the algebraic space $X*_H G$. Consider the $H$-equivariant map $i\colon
X \arr X*_H G$ given by $x\mapsto (x,1)$.  We associate to $F$ the
$H$-torsor $E$ over $S$ defined by the 
pull-back diagram
\[
  \label{e.pb}
\xymatrix{
  E\ar[r]\ar[d]  & X\ar[d]^i\\
        F\ar[r]  & X*_H G.  \\
}
\]

It is not difficult to see that these operations give an equivalence
of categories between $[X/H]$ and $[X*_H G/G]$.
\end{proof} 

Now let $F_{[X/G]}^{\mathrm{spl}}$ be the subfunctor of $F_{[X/G]}$ 
defined as follows. For any field $K/k$, $F_{[X/G]}^{\mathrm{spl}}(K)$ 
consists of diagrams~\eqref{e.functor}, where 
$\pi \colon T  \arr \Spec(K)$ is a split torsor. 

Following~\cite{bf2}, we define the \emph{functor of orbits}
$\Orb_{X, G}$ by $\Orb_{X, G}(K)$ = set of $G(K)$-orbits in $X(K)$. 

\begin{lemma} \label{lem2.2}
The functors $F_{[X/G]}^{\mathrm{spl}}$ and $\Orb_{X, G}$ are isomorphic.
\end{lemma}

\begin{proof} Recall that a torsor $\pi \colon T  \arr \Spec(K)$ is split 
if and only if there exists a section $s \colon \Spec(K)  \arr T$. 

Now we associate the $G(K)$-orbit of the 
$K$-point $\psi s \colon \Spec(K)  \arr X$ to
the object~\eqref{e.functor} of $F_{[X/G]}^{\mathrm{spl}}$.
Note that while the $K$-point $\psi s \colon \Spec(K)  \arr X$
depends on the choice of $s$, its $G(K)$-orbit does not, 
since any other section $s'$ of $\pi$ can be obtained from $s$ by 
translating by an element of $G(K)$. Thus we have defined a map
$\Orb_{X, G}(K)  \arr F_{[X/G]}^{\mathrm{spl}}(K)$ for each $K/k$; it is easy
to see that these maps give rise to a morphism of functors
\begin{equation} \label{e.orb}
\Orb_{X, G}  \arr F_{[X/G]}^{\mathrm{spl}} \, .
\end{equation}

To construct the inverse map, note that a $K$-point  
$p \colon \Spec(K)  \arr X$ of $X$, gives rise to a $G$-equivariant morphism 
$\psi$ from the split torsor $T = G \times \Spec(K)$ to $X$ defined 
by $\psi \colon (g,x)  \arrto g \cdot x$. This morphism represents 
an object in $F_{[X/G]}^{\mathrm{spl}}(K)$.
Translating $p \in X(K)$ by $g \in G(K)$ 
modifies $\psi$ by composing it with an automorphism of $T$ 
given by translation by $g$: 
\[ \xymatrix{ T \ar@{->}[r]^{\times g}  \ar@{->}[rd]^{\pi} & 
T \ar@{->}[r]^{\psi} \ar@{->}[d]^{\pi} &  X \cr & \Spec(K) &  } \]
It is now easy to see that the resulting map 
$F_{[X/G]}^{\mathrm{spl}}  \arr \Orb_{X, G}$ is a morphism of functors,
inverse to~\eqref{e.orb}.
\end{proof}

Recall that a linear algebraic group $G/k$ is called \emph{special}
if every $G$-torsor over $\Spec(K)$ is split, for every 
field $K/k$. 

\begin{corollary} \label{cor3.3} Consider the action of a special 
linear algebraic group $G/k$ on an algebraic space $X$ locally of
finite type $k$. Then

\begin{enumeratea}

\item The functors $F_{[X/G]}$ and $\Orb_{X, G}$ are isomorphic.

\item $\ed{[X/G]} \le \dim X$.

\end{enumeratea}

\end{corollary}

\begin{proof}
(a) Since $G$ is special, $F_{[X/G]}^{\mathrm{spl}} = F_{[X/G]}$. Now apply
Lemma~\ref{lem2.2}.

\smallskip
(b) Let $F_X$ be the functor $K  \arr X(K)$. Then sending a point $p \in X(K)$
to its $G(K)$-orbit induces a surjective morphism of functors 
$F_X  \arr \Orb_{X, G}$. Hence,
   \begin{align*}
   \ed {[X/G]} &= \ed \Orb_{X, G}\\
   &\le \ed F_X\\
   &= \dim(X).\qedhere
   \end{align*}
\end{proof}

\begin{corollary} \label{cor.finite}
Let $G/k$ be a linear algebraic group and let $X/k$ be an algebraic 
space, locally of finite type over $k$ equipped with a $G$-action.
Then $\ed \, [X/G]<\infty$.
\end{corollary}

\begin{proof}
Let $\rho\colon G  \arr \GL_r$ be an embedding and $Y = X *_G \GL_r$.
By Lemma~\ref{lem2.1} the stacks $[X/G]$ and $[Y/\GL_r]$ are 
isomorphic. Since $\GL_r$ is special, Corollary~\ref{cor3.3} 
tells us that $\ed (X/G)=\ed(Y/\GL_r)\leq \dim Y < \infty$.
\end{proof}

Another consequence of Proposition~\ref{p.isotropic} is the following
``classical'' theorem (see~\cite{bf1} for another proof).  

\begin{theorem} 
\label{t.GenFree}
Let $G$ be a linear algebraic group over a field $k$ admitting a
generically free representation on a vector space $V$.  Then
$$
\ed \cB G\leq \dim V-\dim G.
$$
\end{theorem}
\begin{proof}
  Let $U$ denote a dense $G$-stable Zariski open subscheme of $V$ on
  which $G$ acts freely. Then $[U/G]$ is an algebraic space of
  dimension $\dim V-\dim G$ and the map $[U/G] \arr \cB G$ is representable
  and isotropic.
\end{proof}

\subsection*{Finiteness}

The main theorem on finiteness of essential dimension 
is now an easy corollary our study of 
quotient stacks and of a result of A. Kresch.

\begin{theorem}\label{thm.finiteness}
  Let $\cX$ be an algebraic stack of finite type over $k$. If for
  any algebraically closed extension $\Omega$ of $k$ and any
  object $\xi$ of $\cX(\Omega)$ the group scheme
  $\underaut_{\Omega}(\xi) \arr \spec\Omega$ is affine, then
  $\ed(\cX/k) < \infty$.
\end{theorem}

\begin{proof}
By a Theorem of Kresch~\cite[Proposition 3.5.9]{kresch} $\cX$ is covered 
by quotient stacks $[X_i/G_i]$. By Corollary~\ref{cor.finite},
$\ed\cX = \max_i \, \ed {[X_i/G_i]} < \infty$.
\end{proof}

Theorem~\ref{thm.finiteness} does not hold without the assumption that
all the $\underaut_{\Omega}(\xi)$ are affine.  For example, by
Theorem~\ref{thm.curves}, $\ed\cM_{1,0}=+\infty$.   The proof of this will
be given in \S\ref{s.Tate}, and we will also see 
(Theorem~\ref{t.Tate}) that $\ed \cB E=+\infty$ if $E$ is the Tate elliptic
curve over the power series field $\CC\dr{t}$.  

\section{The essential dimension of a smooth \dm stack}
\label{s.generic}

The goal of this section is to prove the following theorem which
allows us, in several of the most interesting cases, to reduce the
calculation of the essential dimension of a stack to that of the
essential dimension of a gerbe over a field.

Recall that if $\cX$ is an algebraic stack over a base scheme $S$,
then $\cX$ is said to be \emph{separated over $S$} when the diagonal
morphism $\Delta: \cX \arr \cX \times_{S} \cX$ is proper. In the
case of a \dm stack, the diagonal morphism is always quasi-finite.  So
the diagonal morphism of a separated \dm stack is finite.

Recall also that the \emph{inertia  stack} $\cI_{\cX} \arr \cX$ is the
fibered product 
$$\cX \times_{\cX \times \cX} \cX$$
mapping to $\cX$ via
the second projection (with  both maps $\cX\arr\cX\times \cX$ given by
the diagonal). The inertia stack  is a group stack, and represents the
functors of isomorphisms of objects:  that is, it is equivalent to the
obvious  fibered category  over $S$  whose objects  are  pairs $(\xi,
\alpha)$, where $\xi$ is an object  over some morphism $T \arr S$, and
$\alpha$ is an automorphism of $\xi$ in $\cX(T)$.

We say that $\cX$ has \emph{finite inertia} when $\cI_{\cX}$ is finite
over $\cX$. A separated \dm stack has finite inertia; however, having finite
inertia is a weaker condition than being separated. For example, when $X$ is a scheme the inertia stack is the identity $X = X$, so $X$ always has finite
inertia, even when it is not separated.

By a result of Keel and Mori (\cite{keel-mori}, see also \cite{conrad-keel-mori}) an algebraic stack locally
of finite type over $\spec k$ with finite inertia has a moduli
algebraic space $\bX$, which is also locally of finite type over
$\spec k$. The morphism $\cX \arr \bX$ is proper.

\begin{theorem}\label{thm:generic}
  Let $k$ be a field of characteristic~$0$, $\cX$ a smooth connected
  \dm stack with finite inertia, locally of finite type over $\spec
  k$. Let $\bX$ the moduli space of $\cX$, $K$ the field of rational
  functions on $\bX$. Denote by $\cX_{K}$ the fibered product $\spec K
  \times_{\bX} \cX$. Then we have
   \[
   \ed(\cX/k) = \dim \bX + \ed(\cX_{K}/K).
   \]
\end{theorem}

\begin{corollary} \label{cor.generic1} 
If $\cX$ is as above and $\cU$ is an open dense substack, 
then $\ed(\cM/k) = \ed(\cU/k)$.
\end{corollary}

\begin{corollary} \label{cor.generic2}
If the conditions of the theorem are satisfied, and the generic object 
of $\cX$ has no non-trivial automorphisms ($\cX$ is an orbifold, 
in the topologists' terminology), then $\ed(\cX/k) = \dim \bX$.
\end{corollary}

\begin{corollary}
\label{cor.curves1}
Assume that $k$ has characteristic~$0$. If $g \geq 3$, or $g = 2$ and
$n \geq 1$, or $g = 1$ and $n \geq 2$, then
   \[
   \ed(\cM_{g,n}/k) = \ed\bigl(\,\overline{\cM}_{g,n}/k\bigr)
   = 3g - 3 + n.
   \]
\end{corollary}

\begin{proof} In all these case the automorphism group of a generic
object of $\cM_{g,n}$ is trivial, so the generic gerbe is trivial, and
$\ed\cM_{g,n} = \dim \cM_{g,n}$. Similarly for $\overline{\cM}_{g,n}$.
\end{proof}

\begin{proof}[Proof of Theorem~\ref{thm:generic}] The equality
$\ed(\cX/k) \geq \dim \bX + \ed(\cX_{K}/K)$ is clear. Let us prove the
opposite inequality.

Let $F$ be an extension of $k$ and $\xi$ an object in $\cX(F)$,
corresponding to a morphism $\xi\colon \spec F \arr \cX$. We need to
show that the essential dimension of $\xi$ is less than or equal to
$\dim \bX + \ed(\cX_{K}/K)$. Of course we can assume that $F$ is
infinite, otherwise $\ed\xi$ would be $0$, in which case we are
done.

We may also assume that $\bX$ is an affine scheme. If it is not so, by
\cite[II, Theorem~6.4]{knutson} the composite $\spec F \xarr{\xi} \cX
\arr \bX$ admits a factorization $\spec F \arr U \arr \bX$, where $U$
is an affine scheme and the morphism $U \arr \bX$ is \'etale. By
substituting $\cX$ with the pullback $\spec U \times_{\bX}\cX$ the
dimension stays the same, while the essential dimension of the generic
gerbe can not increase.

We proceed by induction on the codimension in $\bX$ of the closure of
the image of the composite $\spec F \xarr{\xi} \cX \arr \bX$. If this
codimension is $0$, then $\xi\colon \spec F \arr \cX$ factors though
$\cX_{K}$, in which case the inequality is obvious. So we can assume
that this codimension is positive, that is, the composite $\spec F
\xarr{\xi} \cX \arr \bX$ is not dominant.

\begin{Claim} 
\label{c.Enlarge}
There a morphism $\spec F\ds{t} \arr \cX$, such that its
restriction $\spec F \subseteq \spec F\ds{t} \arr \cX$ is isomorphic
to $\xi$, and such that the image of the composite $\spec F\ds{t}
\arr\cX \arr \bX$ consists of two distinct points.  \end{Claim}

By Schlessinger's theorem, there exists a local complete noetherian
$F$-algebra $A$ with residue field $F$ and is a formal versal
deformation of $\xi$ defined on $A$. Since $\cX$ is not obstructed the
ring $A$ is a power series ring $F\ds{t_{1}, \dots, t_{m}}$. The
composite $\spec A \arr \cX \arr \bX$ is dominant; since $F$ is
infinite, if $a_{1}$, \dots,~$a_{m}$ are general elements of $F$ and
the homomorphism $A \arr F\ds{t}$ is defined by sending $t_{i}$ to
$a_{i}t$, the composite $\spec F\ds{t} \arr \spec A \arr \cX$ has the
required properties.

\begin{Claim}
\label{c.Shrink}
 There exists a complete discrete valuation subring $R
\subseteq F\ds{t}$ and fraction field $L \subseteq F\dr{t}$, such that
the following properties hold:

\begin{enumeratea}

\item $t \in R$,

\item the residue field of $R$ has transcendence degree over $k$ at
most equal to $\dim\bX + \ed(\cX_{K}/K)$, and

\item The composite morphism $\spec F\dr{t} \subseteq \spec F\ds{t}
\arr \cX$ factors through $\spec L$.

\end{enumeratea} \end{Claim}

Consider the composite $\spec F\dr{t} \subseteq \spec F\ds{t} \arr
\cX$. The closure of its image in $\bX$ has a codimension that is less
than the codimension of the closure of $\spec F$; hence by induction
hypothesis there exists an intermediate field $k \subseteq E \subseteq
F\dr{t}$ such that $\spec F\dr{t} \arr \cX$ factors through $\spec E$,
and $\trdeg_{k}L$ is at most $\dim\bX + \ed(\cX_{K}/K)$.

Set $L \eqdef E(t)$ and $R \eqdef L \cap F\ds{t}$; clearly $R$ is a
discrete valuation subring in $L$ containing $t$. We claim that the
transcendence degree of the residue field $R/\mathfrak{m}_{R}$ over
$k$ is less than \[ \trdeg_{k}L \leq \dim\bX + \ed(\cX_{K}/K) + 1.  \]
This is elementary: if $s_{1}$, \dots,~$s_{r}$ are elements of $R$
whose images in $R/\mathfrak{m}_{R}$ are algebraically independent
over $k$, then it is easy to check that $s_{1}$, \dots,~$s_{r}$, $t$
are algebraically independent over $k$. Thus $R$ and $L$ satisfy all
the conditions of the claim, except completeness. By completing $R$ we
prove the claim.

\begin{Claim}
\label{c.Factor}
The morphism $\spec F\ds{t} \arr \cX$ factors through
$\spec R$.  \end{Claim}

This claim implies that $\xi\colon \spec F \arr \cX$ factors through
$\spec(R/\mathfrak{m}_{R})$, which shows that
   \begin{align*}
   \ed\xi
   &\leq \trdeg\bigl((R/\mathfrak{m}_{R})/k\bigr)\\ &\leq \dim\bX +
   \ed(\cX_{K}/K),
   \end{align*}
thus proving the Theorem.

To prove the claim, let us first show that the morphism $\spec F\ds{t}
\arr \bX$ factors through $\spec R$. This is trivial, since $\bX$ is
an affine scheme: if $\bX = \spec A$, the homomorphism $A \arr
F\dr{t}$ corresponding to the composite $\spec F\dr{t} \subseteq \spec
F\ds{t} \arr \cX \arr \bX$ factors through $F\ds{t}$ and also through
$L$, so its image in contained in $L \cap F\ds{t} = R$.

Now denote by $\cX_{R}$ the normalization of the reduced pullback
$(\spec R \times_{\bX} \cX)_{\red}$; by a well known theorem
of Nagata, stating that the normalization of a complete local integral
domain is finite, this is finite over $(\spec R \times_{\bX}
\cX)_{\red}$. The restriction of $\cX_{R}$ to $\spec L
\subseteq \spec R$ coincides with $(\spec L \times_{\bX}
\cX)_{\red}$; hence the morphism $\spec L \arr \cX$ yields a
morphism $\spec L \arr \cX_{R}$. Thus the moduli space $\bX_{R}$ of
$\cX_{R}$, which is integral and finite over $\spec R$, admits a
section over $\spec L$; hence $\bX_{R} = \spec R$.

Since $\spec F\ds{t}$ is normal and the morphism $\spec F\ds{t} \arr
(\spec R \times_{\bX} \cX)_{\red}$ induced by $\spec F \ds{t}
\arr \cX$ is dominant, the morphism $\spec F \ds{t} \arr \cX$ factors
through $\cX_{R}$.

Let $X_{0} \arr \cX_{R}$ be an \'etale map, where $X_{0}$ is a scheme.
Since $R$ is complete, hence henselian, $X_{0}$ contains a component
of the form $\spec R_{0}$, where $R_{0}$ is a discrete valuation ring
which is a finite extension of $R$. This component dominates
$\cX_{R}$, so we can assume $X_{0} = \spec R_{0}$. Set $X_{1} =
X_{0}\times_{\cX_{R}} X_{0}$; we have that $X_{1} = \spec R_{1}$,
where $R_{1}$ is a product of discrete valuation rings. The stack
$\cX_{R}$ has a presentation $X_{1} \double X_{0}$. If we set $\spec
F\ds{t}\times_{\cX_{R}} \spec R_{0} = \spec A$, we have that $A$ is a
product of discrete valuation rings, each of them an \'etale
extension of $F\ds{t}$; this implies that the image of $t \in R$ in
$A$ is a uniformizing parameter in all of the factors of $A$; and this
implies that $t$ is also a uniformizing parameter in $R_{0}$. So
$R_{0}$ is \'etale over $R$, because the characteristic of the base
field is $0$.

Now consider the morphism $\spec L \arr \cX_{R}$: set $\spec L_{0} =
\spec L \times_{\cX_{R}} X_{0}$ and $\spec L_{1} = \spec L
\times_{\cX_{R}} X_{1} = \spec (L_{0}\otimes_{L} L_{0})$. Set also
$S_{i} = R_{i} \otimes_{R} F\ds{t}$ and $M_{i} = L_{i}\otimes_{L}
F\dr{t}$ for $i = 0$ or $1$. The ring $M_{i}$ is a product of fields,
$S_{i}$ a product of discrete valuation rings.

The pullback of the commutative diagram
   \[
   \xymatrix{
   {}\spec L_{1} \ar@<.5ex>[r]\ar@<-.5ex>[r]\ar[d]
      &{}\spec L_{0} \ar[r]\ar[d] &{}\spec L\ar[d]\ar[rd]&\\
   {}\spec R_{1} \ar@<.5ex>[r]\ar@<-.5ex>[r] & {}\spec R_{0} \ar[r]
      &{}\cX_{R} \ar[r] &\spec R
   }
   \]
via the morphism $\spec F\ds{t} \arr \spec R$ is isomorphic to the diagram
   \[
   \xymatrix{
   {}\spec M_{1} \ar@<.5ex>[r]\ar@<-.5ex>[r]\ar[d]
      &{}\spec M_{0} \ar[r]\ar[d] &{}\spec F\dr{t}\ar[d]\ar[rd]&\\
   {}\spec S_{1} \ar@<.5ex>[r]\ar@<-.5ex>[r] & {}\spec S_{0} \ar[r]
      &{}\cX_{F\ds{t}} \ar[r] &\spec F\ds{t}\hsmash{.}
   }
   \]
   But the morphism $\spec F\dr{t} \arr \cX_{F\ds{t}}$ extends to a
   morphism $\spec F\ds{t} \arr \cX_{F\ds{t}}$; and this implies that
   $M_{0}$ is unramified over $F\dr{t}$ with respect to the canonical
   valuation of $F\dr{t}$. Hence $L_{0}$ is unramified over $L$; if we
   denote by $T_{i}$ the normalization of $R_{i}$ in $L_{i}$ we have
   that $T_{0}$ is \'etale over $R$, and $T_{1} = T_{0} \otimes_{R}
   T_{0}$. The diagram
   \[
   \xymatrix{
      {}\spec L_{1} \ar@<.5ex>[r]\ar@<-.5ex>[r]\ar[d]
      &{}\spec L_{0} \ar[d]\\
   {}\spec R_{1} \ar@<.5ex>[r]\ar@<-.5ex>[r] & {}\spec R_{0}
   }
   \]
extends to a diagram
   \[
   \xymatrix{
      {}\spec T_{1} \ar@<.5ex>[r]\ar@<-.5ex>[r]\ar[d]
      &{}\spec T_{0} \ar[d]\\
   {}\spec R_{1} \ar@<.5ex>[r]\ar@<-.5ex>[r] & {}\spec R_{0}
   }
   \]
   which defines the descent data for a morphism $\spec R \arr
   \cX_{R}$ extending $\spec L \arr \cX_{R}$. This proves the theorem.
\end{proof}

\begin{remark}\label{r.generic-gerbe}
The stack $\cX_{K}$ that appears in the statement of Theorem~\ref{thm:generic} can be defined in much greater generality. Let $\cX$ be a locally noetherian integral algebraic stack. It is easy to see that all dominant maps $\spec K \arr \cX$ are equivalent, in the sense of \cite[Definition~5.2]{LMB}, hence they define a point of $\cX$, called the \emph{generic point} of $\cX$. Then the stack $\cX_{K}$ of Theorem~\ref{thm:generic} is the gerbe of $\cX$ at its generic point, in the sense of \cite[\S 11.1]{LMB}. This is naturally called the \emph{generic gerbe} of $\cX$.
\end{remark}

Theorem~\ref{thm:generic} is false in general for algebraic stacks, and also for
singular \dm stacks.

\begin{examples}\hfil

\begin{enumeratea}

\item Let $k$ be any field. Let $G \eqdef \ga \times \ga$ act on
  $\AA^{3}$ be the formula $(s,t)(x,y,z) = (x+sz, y+tz, z)$, and
  define $\cX \eqdef [\AA^{3}/G]$. Let $H \subseteq \AA^{3}$ be the
  hyperplane defined by the equation $z = 0$.  Then $\cX$ is the union
  of the open substack $[(\AA^{3}\setminus H)/G] \simeq \AA^{1}
  \setminus \{0\}$ and the closed substack $[H/G] \simeq \AA^{2}
  \times \cB G$; hence its essential dimension is $2$, its generic
  gerbe is trivial, and its dimension is $1$.

\item Let $r$ and $n$ be integers, $r > 1$. Assume that 
the characteristic of $k$ is prime to $r$. Let 
$X \subseteq \AA^{n}$ be the hypersurface defined by the equation 
$x_{1}^{r} + \dots + x_{n}^{r} = 0$. Let $G \eqdef \mmu_{r}^{n}$ 
act on $X$ via the formula
\[
(s_{1}, \dots, s_{n})(x_{1}, \dots, x_{n}) = (s_{1}x_{1}, \dots, s_{n}x_{n}).
\]
Set $\cX = [X/G]$. Then $\cX$ is the union of $[(X \setminus \{0\})/G]$, 
which is a  quasi-projective scheme of dimension $n-1$, 
and $[\{0\}/G] \simeq \cB \mmu_{r}^{n}$, which has essential dimension $n$.

\item The following example shows that 
Corollary~\ref{cor.generic1} fails even for quotient stacks 
of the form $[X/G]$, where $X$ is a complex affine variety 
and $G$ is a connected complex reductive linear algebraic group. 

Consider the action of $G = \GL_n$ on the affine space $X$ of all
$n \times n$-matrices by multiplication on the left. 
Since $G$ has a dense orbit, and the stabilizer of a non-singular 
matrix in $X$ is trivial, we have
   \[
   \text{$\ed$(generic point of $[X/G]) = 0$}.
   \]
On the other hand, let $Y$ be the locus of matrices of rank~$n-1$, which 
forms a locally closed subscheme of $X$. There is a surjective 
$\GL_n$-equivariant morphism $Y \arr \PP^{n-1}$, sending 
a matrix $A$ into its kernel, 
which induces an isotropic morphism $[Y/G] \arr \PP^{n-1}$. 
Hence by Proposition~\ref{p.observation} we have
   \[
   \ed {[X/G]} \geq \ed {[Y/G]} \geq n-1.
   \]

It is not hard to see that the essential dimension of $[X/G]$ is the maximum of all the dimensions of Grassmannians of $r$ planes in $\CC^{n}$, which is $n^{2}/4$ if n is even, and $(n^{2} - 1)/4$ if $n$ is odd.
\end{enumeratea}
\end{examples}

\begin{question}
One could ask to what class of curves we may apply
Theorem~\ref{thm:generic}. More specifically, 
let $\fM_{g}$ be the stack over $\spec k$, whose objects over 
a $k$-scheme $S$ are flat proper finitely presented maps 
$C \arr S$, whose geometric fibers are reduced locally complete 
intersection irreducible curves of geometric genus~$g$. 
The stack $\fM_{g}$ is an irreducible locally finitely 
presented smooth stack over $\spec k$, but it is not \dm. 
It is not of finite type, either, and it is easy to see 
that $\ed\fM_{g} = \infty$ for all $g$.

However, $\fM_{g}$ will contain the largest open substack 
$\widetilde{\fM}_{g}$ with finite inertia: it follows from 
Theorem~\ref{thm:generic} that $\ed\widetilde{\fM}_{g} = \ed \cM_{g}$. 
Is there a good description of $\widetilde{\fM}_{g}$? 
Does it contain all curves with finite automorphism groups?
\end{question}

\section{Gerbes}
\label{s.gerbes}

In this section, we address the problem of computing the essential
dimension a gerbe over a field $K$. The general problem seems
difficult; however we do have a formula in the case where the gerbe is
banded by $\mu_{p^n}$ for $p$ a prime.

Let $G$ be a sheaf of abelian groups in the category $\Aff_{K}$. A
gerbe $\cX$ over $K$ is said to be \emph{banded by $G$} if for any
affine $K$-scheme $S$ and any object $\xi$ of $\cX$ there is an
isomorphism of group schemes $G_{S} \simeq \underaut_{S}(\xi)$, which
is compatible with pullbacks, in the obvious sense. (Here
$\underaut_{S}(\xi)$ denotes the group scheme of automorphisms of
$\xi$ over $S$.) A gerbe banded by $G$ is neutral if and only if it is equivalent to the classifying stack $\cB_{K} G$.

More generally, \cite{G} contains a notion of gerbe banded by $G$ when
$G$ is not abelian; but we do not need the added generality, which
makes the definition considerably more involved.

There is a natural notion of equivalence of gerbes banded by $G$; the
set of equivalence classes is in natural bijective correspondence with
the group $\H^{2}(K, G)$. The identity is the class of the neutral gerbe
$\cB_{K}G$.

\begin{para} Let $K$ be a field and let $\GG_m$ denote the
  multiplicative group scheme over $K$. Recall that the group
  $\H^2(K,\GG_m)$ is canonically isomorphic to the Brauer group
  $\Br(K)$ of Brauer equivalence classes of central simple algebras
  (CSAs) over $K$.  By Wedderburn's structure theorem, any CSA over
  $K$ isomorphic to the matrix algebra $\rM_n(D)$ for $D$ a division
  algebra over $K$ which is unique up to isomorphism.  Moreover, if
  $A$ and $B$ are two Brauer equivalent CSAs, the division algebras
  $D$ and $E$ corresponding to $A$ and $B$ respectively are
  isomorphic.  For a class $[A]\in\Br(K)$, the \emph{index} of $A$ is
  $\sqrt{\dim_K D}$.

Let $n$ denote a non-negative integer and $\alpha\in \H^2(K,\mu_n)$.
We define the \emph{index} $\ind \alpha$ to be the index of the image on
$\alpha$ under the composition
$$
\H^2(K,\mu_n) \arr \H^2(K,\GG_m)\stackrel{\cong} \larr \Br(K).
$$ 
Note that the index of $\alpha$ is the smallest integer $d$ such that 
$\alpha$ is in the image of the (injective) connecting homomorphism
\begin{equation}
\label{e.partial}
\partial\colon\H^1(K,\PGL_d) \arr \H^2(K,\mu_d)
\end{equation}
arising from the short-exact sequence
\begin{equation}
\label{e.pgl}
1\arr \mu_d \arr \SL_d \arr \PGL_d \arr 1.
\end{equation}

The \emph{exponent} $\ord([A])$ of a class $[A]\in \Br K$
is defined to be its order
in the Brauer group.  Note also that the exponent $\ord([A])$ always
divides the index $\ind([A])$~\cite[Theorem 4.4.5]{herstein}.
\end{para}

The next two sections will be devoted to the proof of the
following result.

\begin{theorem}
\label{t.edGerbe}
  Let $\cX$ be a gerbe over a field $K$ banded by $\mu_{p^m}$ for $p$
  a prime and $k$ a positive integer.   Then 
$$
\ed \cX = \ind {[\cX]}.
$$
\end{theorem}

\section{Canonical dimension of smooth proper varieties}
\label{s.canonical}

The following lemma is well known.

\begin{lemma}\label{lem:compose-rational}
  Let $X$, $Y$ and $Z$ be varieties over $K$. Assume that $Y$ is
  smooth and $Z$ is proper. If there exist rational maps $\alpha \colon
X \darr
  Y$ and $\beta \colon Y \darr Z$, then there exist a rational 
map $\gamma \colon X \darr Z$.
\end{lemma}

\begin{proof} Immediate from Nishimura's lemma; cf.~\cite{nishimura}
or~\cite[Proposition A.6]{ry}.
\end{proof}

\begin{definition}
\label{e.Equivalence}
Two smooth proper varieties $X$ and $Y$ are
\emph{$\rme$-equivalent} (or simply \emph{equivalent})
if there exist rational maps $X \darr Y$ and $Y \darr X$. 
\end{definition}

From Lemma~\ref{lem:compose-rational} above, this is in fact an
equivalence relation. 

\begin{definition}
  If $X$ and $Y$ are smooth proper varieties over $K$, let $\rme(X,Y)$
  denote the least dimension of the closure of the image of a rational
  map $X \darr Y$. We set $\rme(X,Y) = +\infty$ if there are no
  rational maps $X \darr Y$ and define
 $\rme(X) = \rme(X,X)$.
\end{definition}

The integer $\rme(X)$ has been introduced in \cite{km}, and is called
\emph{the canonical dimension of $X$}. (In the case where $X$ is 
a Brauer-Severi variety of dimension $n-1$, this number coincides 
with the canonical dimension of the class of $X$ in $H^1(K, \PGL_n)$,
as defined in~\cite{ber}.) 

\begin{lemma}
  Let $X$, $X'$, $Y$ and $Y'$ be smooth proper varieties over $K$,
  such that $X$ is equivalent to $X'$ and $Y$ is equivalent to $Y'$.
  Then $\rme(X,Y) = \rme(X', Y')$.
\end{lemma}

\begin{proof}
  Let $f\colon X \darr Y$ be a rational map such that $\dim V =
  \rme(X,Y)$, where $V \subseteq Y$ is the closure of the image of
  $f$. From Lemma~\ref{lem:compose-rational}, there exists a rational
  map $X' \darr V$; by composing this with the embedding $V \subseteq
  Y$, we see that there a rational map $X' \darr Y$ whose image has
  dimension at most $\dim V = \rme(X,Y)$. Hence $\rme(X', Y) \leq
  \rme(X, Y)$.

  On the other hand, from the same lemma we see that there a rational
  map $V \arr Y'$; this can be composed with the dominant rational map
  $X \darr V$ to obtain a rational map $X \arr Y'$ whose image has
  dimension at most $\dim V = \rme(X,Y)$. Hence $\rme(X, Y') \leq
  \rme(X, Y)$. We deduce that
   \[
   \rme(X', Y') \leq \rme (X', Y) \leq \rme(X, Y);
   \]
   by symmetry, $\rme(X,Y) = \rme(X', Y')$.
\end{proof}

\begin{corollary}
If $X$ and $Y$ are equivalent smooth proper varieties over $K$, then
   \[
    \rme(X) = \rme(Y) = \rme(X,Y).
   \]
\end{corollary}

When we have resolution of singularities for varieties over $K$, then
$\rme(X)$ can also be defined as the least dimension of a smooth
proper variety in the equivalence class of $X$.

\section{The essential dimension of a gerbe}
\label{s.ed-gerbe}

Let $K$ be a field and let $n$ be an integer with $n > 1$.  Let $\cX
\arr \Spec K$ be a gerbe banded by $\mu_n$ with index $d$. Let $P \arr
\Spec K$ be a Brauer--Severi variety of dimension $d-1$ whose class in
$\H^1(K,\PGL_d)$ maps to $[\cX]\in \H^2(K,\mu_n) \subseteq \Br K$ under the connecting homomorphism $\H^1(K,\PGL_d) \arr \Br K$.

\begin{theorem}
\label{t.eP}
  $\ed\cX = \rme(P)+1$.
\end{theorem}

Before proving the theorem we will prove the following easy
corollaries.

\begin{corollary}
\label{c.EdEquiv}
  If $\cX_1$ and $\cX_2$ are gerbes over $\Spec K$ banded by
  $\mu_{n_1}$ and $\mu_{n_2}$ respectively whose cohomology classes in
  $\Br K$ are the same, then $\ed\cX_1=\ed\cX_2$. 
  \begin{proof}
    Clear.
  \end{proof}
\end{corollary}

\begin{corollary} Theorem~\ref{t.edGerbe} holds.  That is, if $n=p^m$,
  with $m>0$, then $\ed\cX=\ind {[\cX]}$. 
\end{corollary}
\begin{proof}
If the index $d$ is $1$, then $\cX$ is neutral.  Thus $\cX=\cB \mu_n$
with $n>1$.  By~\cite[Example 2.3]{bf1}, $\ed \cB \mu_n=1$.  

Assume then that $d>1$.  Then the class of $\cX$ in  $\Br K$ is also
represented by a gerbe banded by $\mu_d$. By
Corollary~\ref{c.EdEquiv}, we can substitute this gerbe for $\cX$, and
assume $d=n$.  

The class of $\cX$ in $\H^2(K,\mu_d)$ comes from a division algebra of
degree $d$.  If $P$ is the associated Brauer--Severi variety then
$\rme(P)=d-1$ by a theorem of Karpenko~\cite[Theorem 2.1]{karpenko} 
(see also~\cite[\S 7.2]{m1}).  
\end{proof}

\begin{proof}[Proof of Theorem~\ref{t.eP}]

Since the exponent of $\partial P$ divides $n$, there exists an
invertible sheaf $\Lambda$ on $P$ whose degree when base changed to
$P_{\overline{K}}\cong \PP^{d-1}_{\overline{K}}$ is $n$.  
For each $K$ scheme $T$, we denote by $\Lambda_T$ the pullback 
of $\Lambda$ to $P_T\eqdef P\times_K T$.    

The gerbe $\cX$ is equivalent to the gerbe whose sections over a
$K$-scheme $T$ consist of pairs $(L,\lambda)$ where $L$ is a an
invertible sheaf on $P_T$, and $\lambda$ is an isomorphism
of sheaves of $\cO_{P_T}$-modules between $L^{\otimes n}$ and
$\Lambda_T$.  We may, therefore, substitute this gerbe for $\cX$ in
proving Theorem~\ref{t.eP}.

Let $P\dual$ denote the dual Brauer--Severi variety.  Since $P$ splits
over $k(P\dual)$ and $P\dual$ splits over $K(P)$, $P$ and
$P\dual$ are in the same $\rme$-equivalence class.  Thus,
$\rme(P)=\rme(P\dual)$.  Each point $\xi\in P\dual$ gives, by
definition, a hypersurface of degree $1$ in $P_{k(\xi)}$, which we
denote by $H_{\xi}$.

\begin{claim}
\label{c.leq} $\ed\cX\leq \rme(P)+1$.   
\end{claim}
\begin{proof}
  Let $F$ be an extension of $K$, and let $(L,\lambda)$ be a class in
  $\cX(\Spec F)$.  The degree of the pullback of $L$ to
  $P_{\overline{F}}\cong \PP^{d-1}_{\overline{F}}$ is $1$.  So
  $\H^0(P_F,L)$ is an $n$-dimensional vector space over $F$. Choose a
  non-zero section of $L$ and let $H\subset P_F$ denote its divisor.
  Then $H$ gives a morphism $\Spec F \arr P\dual$.  We know that there
  exists a rational map $P\dual\darr P\dual$ whose image has dimension
  $\rme(P\dual)=\rme(P)$; call all $V$ the closure of its image.  Since
  there is a rational map $P\dual\darr V$ and
  $P\dual(F)\neq\emptyset$, we also have $V(F)\neq\emptyset$.  Chose a
  morphism $\Spec F \arr P\dual$ whose image is contained in $V$, and
  call its image $\xi$.  The transcendence degree of $K(\xi)$ over $K$
  is at most $\rme(P)$. Then $k(\xi)\subseteq F$ and the pullback of
  $\cO_{P_{k(\xi)}}(H_{\xi})$ to $P_{F}$ is $\cO_{P_{F}}(H) \simeq L$.
  Fix an isomorphism of invertible sheaves $\mu\colon
  \cO_{P_{k(\xi)}}(H_{\xi})^{\otimes n} \simeq \Lambda_{k(\xi)}$: the
  pullback of $\mu$ gives an isomorphism of invertible sheaves
  $L^{\otimes n} \simeq \Lambda_{F}$, which will differ from $\lambda$
  by an element $a \in F^{*}$. Then $(L,\lambda)$ is clearly defined
  over the field $k(\xi)(a)$, whose transcendence degree over $K$ is
  at most $\rme (P) + 1$.
\end{proof}

Now we prove that $\ed\cX\geq e(P)+1$.

Let $S \subseteq P$ be a divisor in the linear system of $\Lambda$;
when pulled back to $P_{\overline{K}} \simeq
\PP^{d-1}_{\overline{K}}$, the hypersurface $S_{\overline{K}}$ has degree~$n$. If $F$ is an extension of $K$, we can determine an element
of $\cX(\spec F)$ by specifying a hyperplane $H \subseteq P_{F}$ and a
rational function $u \in k(P_{F})$ whose divisor is $S_{F} - nH$
(here, as in what follows, we will write $S_{F}$ to indicate the
pullback of $S$ to $P_{F}$): the line bundle is $\cO_{P_{F}}(H)$, and
the isomorphism
   \[
   \cO_{P_{F}}(H)^{\otimes n} = \cO_{P_{F}}(n H) \simeq \cO_{P_{F}}(S_{F})
   = \Lambda_{F} 
   \]
   is given by multiplication by $u$. Every element of $\cX(\spec F)$
   is isomorphic to one arising in this way. The hyperplane $H$ and
   the rational function $u$ are not unique, but it is easy to see
   that the class of $u$ in $k(P_{F})^{*}/k(P_{F})^{*n}$ is uniquely
   determined by the element of $\cX(\spec F)$. This gives us an
   invariant, which is functorial in $F$.

   Let $P\dual \darr P\dual$ be a rational function whose image has
   dimension $\rme(P)$, and call $V$ the closure of its image in
   $P\dual$. The generic point $\spec k(V) \arr V \subseteq P\dual$
   gives us a rational point $\xi$ of $P^{\vee}_{k(V)}$, corresponding
   to a hyperplane $H_{\xi} \subseteq P_{k(V)}$; let $u$ be a rational
   function on $P_{k(V)}$ whose divisor is $S_{k(V)} - n H_{\xi}$.
   Consider the element $\alpha \in \cX\bigl(\spec k(P)(t)\bigr)$
   determined by the rational function $tu$, whose divisor is
   $S_{k(V)(t)} - n H_{\xi}$. We claim that $\alpha$ can not come from
   an extension of $K$ whose transcendence degree is less than
   $\rme(P) + 1$.

   Let $F$ be a subfield of $k(V)(t)$ containing $K$ such that
   $\alpha$ is defined over $F$; we need to show that the
   transcendence degree of $F$ over $K$ is $\rme(P) + 1$. Let $Z$ be
   an integral variety over $K$ with quotient field $F$. Since $\cX$
   has a section on $F = k(Z)$ there exists a hyperplane $H \subseteq
   P_{k(Z)}$, which gives a rational map $\psi\colon Z \darr
   P^{\vee}$. By blowing up the base locus we can assume that $\psi$
   is a morphism.

   Consider the rational map $\phi\colon V \times_{K} \PP^{1} \darr Z$
   corresponding to the embedding $k(Z) \subseteq k(V)(t)$. The class
   of $tu$ in $\inv{k(P_{k(V)(t)})}$ comes from $\inv{k(P_{k(Z)})}$;
   hence there exist rational functions $v \in k(P_{k(Z)})^{*}$ and $w
   \in k(V_{k(P)(t)})^{*}$ such that
   \[
   tu = w^{n}(\id_{P}\times \phi)^{*}v \in k(P_{k(V)(t)})
      = k(P\times V\times\PP^{1}).
   \]
   Now, the valuation of $tu$ along the divisor $D \eqdef P \times V
   \times\{0\} \subseteq P \times V \times \PP^{1}$ is $1$; since $n >
   1$, the valuation of $(\id_{P}\times \phi)^{*}v$ at $D$ can not be
   $0$. Hence $D$ can not dominate $P \times Z$ under the map $\id_{P}
   \times \phi$, or, equivalently, $V \times \{0\}$ can not dominate
   $Z$ under $\phi$. The restriction of $\phi$ gives a regular
   function $V = V \times \{0\} \arr Z$. The composite
   \[
   V \stackrel{\phi}\longrightarrow Z 
   \stackrel{\psi}\longrightarrow P^{\vee},
   \]
   which is well defined because $\psi$ is a morphism, has image
   closure of dimension less than $\dim Z$: hence, if $\dim Z \leq
   \rme(P) = \rme(P^{\vee})$, by composing with the dominant rational
   map $P\dual \darr V$ we obtain a rational map $P\dual \darr P\dual$
   with image closure less than $\rme(P\dual)$, a contradiction. So
   the dimension of $Z$, which equals the transcendence degree of $F$
   over $K$, is $\rme(P) + 1$, as claimed.

   This completes the proof of Theorem~\ref{t.eP} and, thus, the proof
   of Theorem~\ref{t.edGerbe}.
\end{proof}

When the index of $P$ is not a prime power, the essential dimension of $P$ 
is smaller than the index. In fact, 
let $m$ be the index of $P$, and consider 
the prime decomposition $m = p_{1}^{a_{1}} \dots p_{r}^{a_{r}}$. 
Then the class of $P$ in $\Br K$ is the product of the
classes $\alpha_{1}$, \dots,~$\alpha_{r}$ of indices 
$p_{1}^{a_{1}}$, \dots,~$p_{r}^{a_{r}}$. 
If $P_{1}$, \dots,~$P_{r}$ are Brauer--Severi varieties with classes 
$\alpha_{1}$, \dots, $\alpha_{r}$, then the splitting fields of 
$P$ are exactly the fields that split all of the $P_{i}$; 
hence $P_{1} \times \dots \times P_{r}$ has a point over 
$k(P)$, and $P$ has a point over 
$k(P_{1} \times \dots \times P_{r})$. 
So $P$ and $P_{1} \times \dots \times P_{r}$ are equivalent, and we have
   \begin{align*}
   \rme(P) &= \rme(P_{1} \times \dots \times P_{r})\\
   & \leq p_{1}^{a_{1}} + \dots + p_{r}^{a_{r}} - r.
   \end{align*}
In \cite{ctkm}, Colliot-Th\'el\`ene, Karpenko and Merkurjev conjecture 
that equality always holds and prove it for $m = 6$.
This can be reformulated in the language of essential dimension 
of gerbes in the following fashion.

\begin{conjecture}\label{conj:ed-product}
If $\cX$ is a gerbe banded by $\mu_{n}$ over a field $K$, 
let $p_{1}^{a_{1}} \dots p_{r}^{a_{r}}$ be the decomposition 
into prime factors of the index of the class of $\cX$ in 
the Brauer group of $K$. Then
   \[
   \ed\cX = p_{1}^{a_{1}} + \dots + p_{r}^{a_{r}} - r + 1.
   \]
\end{conjecture}

When the index is $6$ this follows from \cite[Theorem 1.3]{ctkm}.

In view of the fact that this holds for $r = 1$, the conjecture 
can be rephrased as follows: if $m$ and $n$ are relatively prime 
positive integers, $\cX$ and $\cY$ are gerbes banded by 
$\mu_{m}$ and $\mu_{n}$, then
   \[
   \ed(\cX \times \cY) = \ed\cX + \ed\cY - 1.
   \]
Or, back to the language of canonical dimension, one could ask the following more general question. Let $X$ and $Y$ be smooth projective varieties over 
a field $K$. Assume that there are no rational functions $X \darr Y$ or $Y \darr X$. Then is it true that $\rme(X \times Y) = \rme(X) + \rme(Y)$? A positive answer to this question would imply the conjecture above.

\section{The essential dimension of $\cM_{g,n}$ for $(g,n) \neq (1,0)$}
\label{s.ed-Mgn}

In this section we complete the proof of Theorem~\ref{thm.curves} when
$(g,n) \neq (1,0)$. By Corollary~\ref{cor.curves1}, it will suffice to
compute $\ed\cM_{0,0}$, $\ed\cM_{0,1}$ and $\ed\cM_{0, 2}$,  
$\ed\cM_{1,1}$ and $\ed\cM_{2,0}$.  

It is easy to see that $\ed \, \cM_{0,1} = \ed \, \cM_{0,2} = 0$.
Indeed, a smooth curve $C$ of genus $0$ with one or two rational points 
over an extension $K$ of $k$ is isomorphic to $(\PP^{1}_{k}, 0)$ or
$(\PP^{1}_{k}, 0, \infty)$, hence it is defined over $k$. 
Alternatively, $\cM_{0,2} = \cB\gm$ and $\cM_{0,1} = 
\cB(\gm \ltimes \GG_{\rma})$, 
and the groups $\gm$ and $\gm \ltimes \GG_{\rma}$ are special 
(and hence have essential dimension $0$).

We will now consider the remaining cases, starting with $\cM_{0, 0}$.

Since $\cM_{0,0}\simeq \cB \PGL_2$, the fact that $\ed\cM_{0,0} = 2$ is
classical.  We recall the following result.

\begin{proposition} Let $k$ be an algebraically closed field of
  characteristic $0$.  Then $\ed\PGL_n=2$ for
$n=2, 3$ and $6$.
\end{proposition}

This is~\cite[Lemma 9.4 (c)]{reichstein}.  However, note that the
proof does not really require the field $k$ to be
algebraically closed of characteristic $0$: it goes through whenever
$\chr k$ does not divide $n$ and $k$ contains a primitive $n$-th root of
unity. Thus we have the following.

\begin{corollary}
  Let $k$ be a field of characteristic not equal to $2$.   Then
  $\ed\cM_{0,0}=\ed\cB\PGL_2=2$.
\end{corollary}

This can also be proved directly very easily: 
the inequality $\ed\cM_{0,0} \leq 2$ holds because 
every smooth curve of genus~$0$ over a field $K$ is a conic 
in $\PP^{2}_{K}$, and can be defined by an equation 
of the type $ax^{2} + by^{2} + x^{2} = 0$ for some 
$a$, $b \in K$, hence is defined over $k(a,b)$. 
The opposite inequality follows from Tsen's theorem.

We will now proceed to compute the essential dimension of $\cM_{1, 1}$.

\begin{proposition} \label{t.elliptic}
If $k$ is a field of characteristic not equal to $2$ or $3$, then 
   \[
   \ed\cM_{1,1}=2.
   \]
\end{proposition}

\begin{proof}
  Every elliptic curves over a field $K$ can be written as a cubic in
  $\PP^{2}_{K}$ with equation $yz = x^{3} + a xz^{2} + bz^{3}$, so it
  is defined over $k(a,b)$. Hence $\ed \cM_{1,1} \leq 2$.

  Let $\cM_{1,1} \arr \AA^1_{k}$ denote the map given by the
  $j$-invariant and let $\cX$ denote the pull-back of $\cM_{1,1}$ to
  the generic point $\Spec k(j)$ of $\AA^1$. Then $\cX$ is banded by
  $\mu_2$ and neutral by~\cite[Proposition 1.4 (c)]{Silverman}, and so
  $\ed \cX = \ed \cB_{k(j)} \mu_{2} = 1$. This implies what we want.
\end{proof}

Ir remains to compute the essential dimension of $\cM_{2, 0}$.
The equality $\ed \cM_{2,0} = 5$ is a special case of 
the following more general result.

\begin{theorem}
\label{t.hyperelliptic}
Let $\cH_g$ denote the stack of hyperelliptic curves of genus $g>1$
over a field $k$ of characteristic $0$ and let $\overline{\cH}_g$
denote its closure in $\overline{\cM}_g$.  Then
$$
\ed \cH_g= \ed \overline{\cH}_g=
\begin{cases}
 2g   & \text{if $g$ is odd,}\\
 2g+1 & \text{if $g$ is even.}\\
\end{cases}
$$
\end{theorem} 

Since $\cH_{2} = \cM_{2,0}$ and $\overline{\cH}_{g} =
\overline{\cM}_{2}$, Theorem~\ref{t.hyperelliptic}
completes the proof of every case of
Theorem~\ref{thm.curves}, except for the fact that $\ed\cM_{1,0}=+\infty$.

\begin{proof} The closure $\overline{\cH}_{g}$ is well known to be
smooth, so by Corollary~\ref{cor.generic1} it is enough to prove the
statement about $\cH_{g}$. Denote by $\bH_{g}$ the moduli space of
$\cH_{g}$; the dimension of $\cH_{g}$ is $2g-1$. Let $K$ be the field
of rational functions on $\bH_{g}$, and denote by $(\cH_{g})_{K}
\eqdef \spec K \times_{\bH_{g}} \cH_{g}$ the generic gerbe of
$\cH_{g}$. From Theorem~\ref{thm:generic} we have \[ \ed \cH_{g} = 2g
- 1 + \ed\bigl((\cH_{g})_{K}/K\bigr), \] so we need to show that
$\ed\bigl((\cH_{g})_{K}/K\bigr)$ is $1$ if $g$ is odd, $2$ if $g$ is
even.

For this we need some standard facts about stacks of hyperelliptic
curves, which we recall.

Call $\cD_{g}$ the stack over $K$ whose object over a $K$-scheme $S$
are pairs $(P\to S, \Delta)$, where $P \to S$ is a conic bundle (that
is, a Brauer--Severi scheme of relative dimension~$1$), and
$\Delta\subseteq P$ is a Cartier divisor which is \'etale of
degree~$2g+2$ over $S$. Every family $\pi\colon C \to S$ in $\cH(S)$
comes with a unique flat morphism $C \arr P$ of degree~$2$, where $P
\arr S$ is a smooth conic bundle; denote by $\Delta \subseteq P$ its
ramification locus. Sending $\pi\colon C \to S$ to $(P \arr S,
\Delta)$ gives a morphism $\cH_{g} \arr \cD_{g}$. Recall the usual
description of ramified double covers: if we split $\pi_{*}\cO_{C}$ as
$\cO_{P} \oplus L$, where $L$ is the part of trace~$0$, then
multiplication yields an isomorphism $L^{\otimes 2} \simeq
\cO_{P}(-\Delta)$. Conversely, given an object $(P \to S, \Delta)$ of
$\cD_{g}(S)$ and a line bundle $L$ on $P$, with an isomorphism
$L^{\otimes 2} \simeq \cO_{P}(-\Delta)$, the direct sum $\cO_{P}
\oplus L$ has an algebra structure, whose relative spectrum is a
smooth curve $C\to S$ with a flat map $C \arr P$ of degree~$2$.

The morphism $\cH_{g} \arr \bH_{g}$ factors through $\cD_{g}$, and the
morphism $\cD_{g} \arr \bH_{g}$ is an isomorphism over the non-empty
locus of divisors on a curve of genus~$0$ with no non-trivial
automorphisms (this is non-empty because $g \geq 2$, hence $2g+2 \geq
5$). Call $(P \arr \spec K, \Delta)$ object of $\cD_{g}(\spec K)$
corresponding the generic point $\spec K \arr \bH_{g}$. It is well
known, and easy to show, that $P(K) = \emptyset$. By the description
above, the gerbe $(\cH_{g})_{K}$ is the stack of square roots of
$\cO_{P}(-\Delta)$, which is banded by $\mu_{2}$. When $g$ is odd then
there exists a line bundle of degree $g+1$ on $P$, whose square is
isomorphic to $\cO_{P}(-\Delta)$; this gives a section of
$(\cH_{g})_{K}$, which is therefore isomorphic to $\cB_{K}\mu_{2}$,
whose essential dimension over $\mu_{2}$ is $1$. If $g$ is even then
such a section does not exist, and the stack is isomorphic to the
stack of square roots of $\omega_{P/K}$, whose class in $\H^{2}(K,
\mu_{2})$ represents the image in $\H^{2}(K, \mu_{2})$ of the class
$[P]$ in $\H^{1}(K, \PGL_{2})$ under the non-abelian boundary map
$\H^{1}(K, \PGL_{2}) \arr \H^{2}(K, \mu_{2})$. According to
Theorem~\ref{t.edGerbe} its essential dimension is the index of $[P]$,
which equals $2$.  \end{proof}

\section{Central Extensions}
\label{s.extensions}

The rest of this paper will rely on our analysis of the following
situation which we recall from the introduction~\eqref{e.extensions}.

Let 
\begin{equation}
\label{e.extensionsb}
1\arr Z \arr G \arr Q\arr 1
\end{equation}
denote an extension of group schemes over a field $k$ with $Z$ central
and isomorphic to $\mu_n$ for some integer $n>1$. 
As in the introduction, 
we define $\ind(G,Z)$ as the maximal value of 
$\ind\bigl(\partial_K(t)\bigr)$ as 
$K$ ranges over all field extensions of $k$ and $t$ ranges over all
torsors in $\H^1(K,Q)$. 

We are now going to prove Theorem~\ref{thm2} from the introduction
which we restate for the convenience of the reader.

\begin{theorem}
\label{thm2b} 
Let $G$ be an extension as in~\eqref{e.extensionsb}.  Assume that $n$
is a prime power. Then
\[
\ed (\cB G/k) \geq \ind(G,Z) - \dim Q.
\]
\end{theorem}

\begin{proof}
Let $K/k$ be a field extension and let $t:\Spec K\to\cB Q$ be a
$Q$-torsor over $\Spec K$. The dimension of $\cB Q$ at the
point $t$ is $-\dim Q$.   Let $\cX$ denote the pull-back in the
following diagram.
$$\xymatrix{
\cX\ar[r]\ar[d] & \Spec K\ar[d]^{t}\\
\cB G\ar[r]     & \cB Q        \\
}
$$
By Proposition~\ref{p.RelDim2}, $\ed(\cX/K)\leq \ed(\cB G/k) + \dim Q$.
On the other hand, since $\cB G$ is a gerbe banded by $Z$ over $\cB
Q$, $\cX$ is a gerbe banded by $Z$ over $\Spec K$.   Therefore, by
Theorem~\ref{t.edGerbe}, $\ed (\cX/K)=\ind\partial_K(t)$.  
By substitution, $\ind\partial_K(t)-\dim Q\leq \ed(\cB G/k)$.   Since this
inequality holds for all field extensions $K/k$ and all $Q$-torsors
$t$ over $K$, the result follows.
\end{proof}

\begin{remark} \label{rem.ctkm}
An affirmative answer to Conjecture~\ref{conj:ed-product} would yield
an inequality similar to the one in~\ref{thm2b} without the assumption
that $n$ is a prime power:  Let $\ind(G,Z)=\prod  p_i^{a_i}$ be the
prime factorization of $\ind(G,Z)$.    The
conjecture would imply that 
\[
\ed (\cB G/k) \ge 1 - \dim Q + \sum (p_i^{a_i}-1).
\]
As remarked in \S 7, the conjecture is a theorem in the
case that the index is $6$ (\cite[Theorem Theorem 1.3]{ctkm}).  We
therefore have that 
\[ \ed (\cB G/k) \ge 4 - \dim Q \, . \]
\end{remark}

\begin{remark} Suppose $G$ is a simple algebraic group whose center $Z$ is
cyclic.  It is tempting to apply Theorem~\ref{thm2b} to the natural
sequence
\[ 1 \arr Z \arr G \arr G^{\mathrm{ad}} \arr 1 \] where the adjoint group
$G^{\mathrm{ad}}$ is $G/Z$. Given a torsor $t \in H^1(K, G^{\mathrm{ad}})$, the central
simple algebra representing $\partial_K(t) \in H^2(K, Z)$ is called {\em
  the Tits algebra} of $t$. The possible values of the index of the
Tits algebra were
studied in~\cite{tits}, where it is denoted by $b(X)$ (for
group of type $X$) and its possible values are listed on p. 1133. A
quick look at this table reveals that for most types these indices are
smaller than $\dim(G)$, so that the bound of Theorem~\ref{thm2b}
becomes vacuous. The only exception are groups of types $B$ and $D$,
in which case Theorem~\ref{thm2b} does indeed, give interesting
bounds; cf. Remark~\ref{rem.clifford}.
\end{remark}

\section{Tate curves and the essential dimension of $\cM_{1,0}$}
\label{s.Tate}

Our first application of Theorem~\ref{thm2b} is to finish the
proof of Theorem~\ref{thm.curves} from the introduction by showing
that $\ed\cM_{1,0}=+\infty$.  

Note that by $\cM_{1,0}$ we mean the moduli stack of genus $1$ curves,
not the moduli stack $\cM_{1,1}$ of elliptic curves (which is
Deligne-Mumford). The objects of $\cM_{1,0}$ are torsors for elliptic
curves as opposed to the elliptic curves which appear as the objects
of $\cM_{1,1}$.  We will now see that these torsors are what causes
the essential dimension to be infinite.

\begin{para} 
\label{p.Tate} Let $R$ be a complete discrete valuation ring with
  function field $K$ and uniformizing parameter $q$.  For simplicity,
  we will assume that $\chr K=0$.  Let $E=E_q/K$ denote
  the Tate curve over $K$~\cite[\S 4]{Silverman}.  This is an elliptic
  curve over $K$ with the property that, for every finite field
  extension $L/K$, $E(L)\cong L^*/q^{\ZZ}$.  It follows that the 
  kernel  $E[n]$ of multiplication by an integer $n>0$ fits
  canonically into a short exact sequence
\[
\label{t.ses}
0 \arr \mu_n \arr E[n] \arr \ZZ/n \arr 0.
\]
Let $\partial\colon\H^0(K,\ZZ/n) \arr \H^1(K,\mu_n)$ denote the connecting
homomorphism.  Then it is well-known (and easy to see) that
$\partial(1)=q\in \H^1(K,\mu_n)\cong K^*/(K^*)^n$.  
\end{para}

\begin{theorem}
\label{t.Tate} Let $E=E_q/K$ denote the Tate curve over a field $K$ as
in ~\eqref{p.Tate}.  Then  
$$
\ed E=+\infty.
$$
\end{theorem}

Theorem~\ref{t.Tate} is an immediate consequence of the following statement.

\begin{lemma}
\label{l.TateTorsion}
  Let $E=E_q$ be a Tate curve as in~\eqref{p.Tate} and let $l$ be a
  prime integer not equal to $\chr R/q$.  Then, for any integer $n>0$, 
  \begin{equation*}
    \ed E[l^n] = l^n.
  \end{equation*}
\end{lemma}
\begin{proof}
We first show that $\ed E[l^n]\geq l^n$.

Let $R'\eqdef R[1^{1/l^n}]$ with fraction field $K'=K[1^{1/l^n}]$.  Since
$l$ is prime to the residue characteristic, $R'$ is a complete
discrete valuation ring, and the Tate curve $E_q/K'$ is the pullback
to $K'$ of $E_q/K$.  Since $\ed (E_q/K')\leq \ed (E_q/K)$, it suffices
to prove the lemma with $K'$ replacing $K$.  In other words, it
suffices to prove the lemma under the assumption that $K$ contains the
$l^n$-th roots of unity.

In that case, we can pick a primitive $l^n$-th root of unity $\zeta$
and write $\mu_{l^n}=\ZZ/l^n$.   Let $L=K(t)$ and consider the 
class $(t)\in \H^1(L,\mu_{l^n})=L^*/(L^*)^n$.  

It is not difficult to see that 
$$
\partial_K(t)=q\cup (t).
$$
It is also not difficult to see that the order of $q\cup (t)$ is
$l^n$ (as the map $\alpha\mapsto \alpha\cup (t)$ is injective by
cohomological purity).
Therefore $\ind(q\cup (t))=l^n$.  It follows that
$\ind(E[l^n],\mu_{l^n})\geq l^n$.  Then, since $\dim\ZZ/l^n=0$,
Theorem~\ref{thm2b} implies that $\ed \cB E[l^n]\geq l^n$. 

To see that $\ed \cB E[l^n]\leq l^n$, note that $E[l^n]$ admits an 
$l^n$-dimensional generically free representation $V=\Ind_{\mu_{l^n}}^{E[l^n]}
\chi$ where $\chi\colon\mu_{l^n} \arr \GG_m$ is the tautological character.
Thus, by Theorem~\ref{t.GenFree}, we have the desired inequality.
\end{proof}

\begin{proof}[Proof of Theorem~\ref{t.Tate}]
  For each prime power $l^n$, the morphism $\cB E[l^n] \arr \cB E$ is
  representable of fiber dimension $1$.  We therefore have
\begin{align*}
\ed E&\geq \ed \cB E[l^n] \\
     &=l^n-1
\end{align*}
for all $n$.
\end{proof}

G.~Pappas pointed out the following corollary.

\begin{corollary}
  Let $E$ be a curve over a number field $K$.  Assume that there is at
  least one prime $\frp$ of $K$ where $E$ has semistable bad reduction.
  Then $\ed E=+\infty$.
\end{corollary}
 
It seems reasonable to make the following guess.

\begin{conjecture}
  If $E$ is an elliptic curve over a number field, then $\ed E=+\infty$.
\end{conjecture}

\begin{remark}
  Note, however, that, if $A$ is a $d$-dimension complex abelian variety, then 
$\ed A=2d$; see ~\cite{brosnan}.
\end{remark}

Now we can complete the proof of Theorem~\ref{thm.curves}.

\begin{theorem}
  \label{t.m1}  Let $k$ be a field. Then $\ed(\cM_{1,0}/k)=+\infty$.   
\end{theorem}
\begin{proof}
  Consider the morphism $\cM_{1,0} \arr \cM_{1,1}$ which sends a genus $1$
  curve to its Jacobian.  Let $F=k\dr{t}$ and let $E$ denote the Tate
  elliptic curve over $F$, which is classified by a morphism $\Spec
  F \arr \cM_{1,1}$.  We have a Cartesian diagram:
  \begin{equation}
    \label{eq:pbe}
\xymatrix{
    \cB E\ar[r]\ar[d] & \cM_{1,0}\ar[d]\\
    \Spec F\ar[r]  & \cM_{1,1}.
}
  \end{equation}
It follows that $+\infty=\ed \cB E\leq \ed\cM_{1,0}$.
\end{proof}

\section{Essential dimension of $p$-groups I}
\label{s.ed-pgroups}

The goal of this section is to prove the following theorem from the
introduction (Theorem~\ref{t.p-groups}).

\begin{theorem}\label{t.p-groups1}
Let $G$ be a $p$-group whose commutator $[G, G]$ is central and cyclic.
Then

\begin{enumeratea}

\item We have
   \[
   \ed_k G \ge \sqrt{|G/\rC(G)|} + \rank \, \rC(G) - 1
   \]
for any base field $k$ of characteristic $\ne p$.

\item Moreover, assume that $k$ contains a primitive root of unity of 
degree $\exp(G)$. Then $G$ has a faithful representation of degree 
   \[
   \sqrt{|G/\rC(G)|} + \rank \, \rC(G) - 1\,.
   \]

\end{enumeratea}
\end{theorem}

Theorem~\ref{t.p-groups} is an immediate consequence of this result,
since part (b) implies
$ \ed_k G \le \sqrt{|G/\rC(G)|} + \rank \, \rC(G) - 1$;
cf. e.g., Theorem~\ref{t.GenFree}.
Our proof of part (a) will rely on the following lemma.

\begin{lemma} \label{lem5.3} Let $G$ be a finite group and $H$ be 
a central cyclic subgroup. Assume that there exists a character 
$\chi \colon G \to k^*$ whose restriction to $H$ is faithful.
Then
\begin{enumerate}

\item $\ed(G) \ge \ed(G/H)$.

\item Moreover, if $H$ has prime order and is not properly contained
in another central cyclic subgroup of $G$ then $\ed(G) = \ed(H) + 1$.

\end{enumerate}
\end{lemma}

\begin{proof} Part (2) is proved in~\cite[Theorem~5.3]{bur} 
in characteristic zero and in~\cite[Theorem 4.5]{kang} 
in prime characteristic. 

To prove (1), let $\phi \colon G \to G/H \hookrightarrow \GL(V)$ 
be a faithful representation of $G/H$. Then 
$\phi \oplus \chi \colon G \to \GL(V \times k)$ is 
a faithful representation of $G$.  Denote the class of the
$G$-action on $V \times k$ by 
$\alpha \in H^1(K, G)$, where $K = k(V \oplus k)^G$. 
Let $\beta$ be the image of $\alpha$ in $H^1(K, G/H)$.  Then 
$\beta$ is given by the induced action of $G/H$ on 
\[ (V \times k)/H \simeq V \times k \, . \] 
Here the quotient map $V \times k \to V \times k$ is given by
$(v, x) \to (v, x^d)$, where $d = |H|$. This shows that induced action of
$G/H$ on $(V \times k)/H$ is again linear. Hence, $\beta$ is a versal
$G/H$-torsor; cf.~\cite[Example 5.4]{gms} 
or \cite[Theorem 3.1]{bur}. We conclude that
\[ \ed(G) = \ed(\alpha) \ge \ed(\beta) = \ed(G/H) \, , \]
as claimed.
\end{proof}

We now proceed with the proof of Theorem~\ref{t.p-groups1}.
We begin with the following reduction.

\begin{lemma} \label{lem.reduction}
In the course of proving Theorem~\ref{t.p-groups1}
we may assume without loss of generality that the center $\rC(G)$ is cyclic.
\end{lemma}

\begin{proof}
Let $Z$ be a maximal cyclic subgroup of $\rC(G)$ containing $[G, G]$.
Then $\rC(G) = Z \oplus W$ for some central subgroup $W$ of $G$. 
Note that $\rank(W) = \rank \, \rC(Z) - 1$.  Moreover, 
$W$ projects isomorphically into the abelian subgroup
$G/Z$. In particular, if $H$ is a cyclic subgroup of $W$ then
after composing this projection with a suitable character of $G/Z$,
we obtain a character $\chi_H \colon G \to k^*$, which is faithful on $H$.
This means that Lemma~\ref{lem5.3} can be used to compare the essential 
dimensions of $G$ and $G/H$.

We will now argue by induction on $|W|$. If $|W| = 1$, 
we are done. For the induction step we will choose
a cyclic subgroup $H \subset W$ (in a way, to be specified below)
and assume that parts (a) and (b) of Theorem~\ref{t.p-groups1} hold
for $\overline{G} = G/H$. Our goal will then be to prove that they 
also hold for $G$.

It is easy to see that $g_1$ and $g_2$ commute in $G$ if and 
only if their images $\overline{g_1}$ and $\overline{g_2}$ commute 
in $\overline{G}$. 
In particular, $\rC(\overline{G}) = \rC(G)/H \simeq Z \oplus W/H$. 
(Here $Z$ projects isomorphically to a cyclic subgroup $\overline{Z}$ of
$\overline{G} = G/H$, and we are identifying $Z$ with $\overline{Z}$).
Consequently,
\begin{equation} \label{e.centralizer}
\text{$|\overline{G}/ \rC(\overline{G})| = |G/\rC(G)|$ and 
$ \rank \, \rC(\overline{G}) = \rank(W/H) + 1$.}
\end{equation}

\smallskip
(a) We choose $H$ to be a subgroup of prime order in $W$. By
induction assumption,
\[ \ed(\overline{G}) \ge \sqrt{|\overline{G}/\rC(\overline{G})|} + 
\rank(\rC(\overline{G})) - 1 \, . \]  
We will now consider two cases.

\smallskip {\bf Case 1.} 
$H$ is properly contained in another cyclic subgroup 
of $W$. In this case $\rank(W/H) = \rank(W)$ and by~\eqref{e.centralizer}
\[ \rank \, \rC(\overline{G}) = \rank(W/H) + 1 = \rank(W) + 1 = \rank \, C(G) \, . \] 
By Lemma~\ref{lem5.3}(1), $ \ed(G) \ge \ed(\overline{G})$ and thus
   \begin{align*}
   \ed(G) &\ge \sqrt{|\overline{G}/\rC(\overline{G})|} + 
   \rank(\rC(\overline{G})) - 1\\
   &= \sqrt{|G|/|\rC(G)|} + \rank \, \rC(G) - 1 \, , 
   \end{align*}as desired.

\smallskip {\bf Case 2.} $H$ is not properly contained in any cyclic subgroup 
of $W$. In this case $\rank(W/H) = \rank(W) - 1$ and by~\eqref{e.centralizer}
\[ \rank \, \rC(\overline{G}) = \rank(W/H) + 1 = \rank(W) = 
\rank \, \rC(G) - 1 \, . \] 
By Lemma~\ref{lem5.3}(2), $ \ed(G) = \ed(\overline{G}) + 1$ and thus
   \begin{align*}
   \ed(G) &= \ed(\overline{G}) + 1\\
   &= \sqrt{|\overline{G}/ \rC(\overline{G})|} + \rank(\rC(\overline{G}))\\
   &= \sqrt{|G|/|\rC(G)|} + (\rank(\rC(G)) -1 ) \, .
   \end{align*}

\smallskip
(b) Here we choose $H$ to be a maximal central cyclic subgroup of $W$
(not necessarily of prime order). By our induction assumption,
$G/H$ has a representation $\rho \colon G/H \hookrightarrow \GL(V)$
of dimension 
   \[
   \sqrt{|\overline{G}/\rC(\overline{G})|} + 
   \rank(\rC(\overline{G})) - 1 \, .
   \]  
Then $\rho \oplus \chi_H$ is a faithful representation of $G$ of dimension
\[ \sqrt{|\overline{G}/\rC(\overline{G})|} + \rank \, \rC(\overline{G}) =
 \sqrt{|G|/|\rC(G)|} + \rank(\rC(G)) - 1 \, ; \]  
cf.~\eqref{e.centralizer}. This shows that Theorem~\ref{t.p-groups1} 
holds for $G$.
\end{proof}

\begin{proof}[Proof of Theorem~\ref{t.p-groups1}(a)]
By Lemma~\ref{lem.reduction} we may assume that $\rC(G)$ is cyclic. 
Set $Z = \rC(G)$. Since we are assuming that
$[G, G] \subset Z$, the quotient $A \eqdef G/Z$ is abelian. 

In this case Theorem~\ref{t.p-groups} reduces to $\ed(G) \ge \sqrt{|A|}$.
By Theorem~\ref{thm2b} it suffices to show that 
\[ \ind(G,Z) \geq \sqrt{|A|} \, . \]
so we will now direct our attention towards computing $\ind(G, Z)$.
 
We will use additive notation for the groups $Z$ and $A$, 
multiplicative for $G$. 
In this situation we can define a skew-symmetric
bilinear form $\omega \colon A \times A  \arr Z$ by
   \[
   \omega(a_1, a_2) = g_1g_2 g_1^{-1} g_2^{-1} \, ,
   \]
where $a_i = g_i$, modulo $Z$, for $i = 1, 2$. (Note that
$\omega(a_1, a_2)$ is independent of the choice of $g_1$ and $g_2$.)
Clearly $g$ lies in $\rC(G)$ if and only if its image $a$ lies in
the kernel of $\omega$, i.e., $\omega(a,b) = 0$ for every $b \in A$. 
Since we are assuming that $\rC(G) = Z$, we conclude that the kernel of
$\omega$ is trivial, i.e., $\omega$ is a {\em symplectic form} on $A$.
It is well known (see for example 
\cite[\S 3.1]{tignol-amitsur}) that the order of $A$, 
which equals the order of $G/\rC(G)$, is a complete square.

Fix a generator $z$ of $Z$. We recall the basic 
result on the structure of a symplectic form $\omega$ on
a finite abelian group $A$ (the proof is 
easy; it can be found, e.g., in~\cite[\S 3.1]{wall}
or~\cite[\S 7.1]{tignol-amitsur}). 
There exist elements $a_1$, \dots,~$a_{2r}$ in $A$ 
and positive integers $d_{1}$, \dots,~$d_{r}$ with 
the following properties.

\begin{enumeratea}
\item $d_{i}$ divides $d_{i-1}$ for each $i = 2$, \dots,~$r$, 
and $d_{r} > 1$.

\item Let $i$ be an integer between $1$ and $r$. 
If $A_{i}$ denotes the subgroup of $A$ generated by $a_i$ and $a_{r+i}$, 
then there exists an isomorphism $\overline{A}_{i} \simeq 
(\ZZ/d_{i}\ZZ)^{2}$ such that 
$a_{i}$ corresponds to (1,0) and $a_{r+i}$ to $(0,1)$.

\item The subgroups $A_{i}$ are pairwise orthogonal 
with respect to $\omega$.

\item $\omega(a_{i}, a_{r+i}) = z^{n/d_{i}}\in Z$.

\item $A = A_{1} \oplus \dots \oplus A_{r}$.
\end{enumeratea}

Then the order of $A$ is $d_{1}^{2} \dots d_{r}^{2}$, 
hence $\sqrt{|A|} = d_{1} \dots d_{r}$.

Let $G_{i}$ be the inverse image of $A_{i}$ in $G$; 
note that $G_{i}$ commutes with $G_{j}$ for any $i \ne j$.

Let $u_{1}$, \dots,~$u_{2r}$ be indeterminates, and 
set $K \eqdef k(u_{1}, \dots, u_{2r})$. 
Identify $Z$ with $\mu_{n}$ by sending $z$ into $\zeta_{n}$. 
Consider the boundary map
   \[
   \partial_{i}\colon \H^{1}(K, A_{i}) \arr \H^{2}(K, Z)
   \]
obtained from the exact sequence
   \[
   1 \arr Z \arr G_{i} \arr A_{i} \arr 1.
   \]

\begin{claim}
There exists a class $\xi_{i} \in \H^{1}(K,A_{i})$ 
such that $\partial_{i}\xi_{i}$ is the class of the cyclic 
algebra $(u_{i}, u_{r+i})_{d_{i}}$ in $\Br K$.
\end{claim}

To see that Theorem~\ref{t.p-groups1}(a) follows from the claim,
consider the commutative diagram
   \[
   \xymatrix{
   1 \ar[r] & Z ^{r} \ar[r]\ar[d]^{m}& {}\prod_{i} G_{i}\ar[r] \ar[d]
     & \prod_{i}A_{i} \ar[r]\ar[d] & 1\\
   1 \ar[r] & Z \ar[r] & G \ar[r] & A \ar[r] & 1
   }
   \]
in which $m$ is defined by the formula 
$m(z_{1}, \dots, z_{r}) = z_{1} \dots z_{r}$, and 
the homomorphism $\prod_{i} G_{i} \arr G$ is induced
by the inclusions $G_{i} \subseteq G$. This yields 
a commutative diagram
   \[
   \xymatrix{
   {}\prod_{i}\H^{1}(K,A_{i}) \ar[r]^{\prod_{i}\partial_{i}} \ar[d]\ar[r]&
   {}\H^{2}(K, Z)^{r}  \ar[d]^{m_{*}}\\
   {}\H^{1}(K,A) \ar[r]^{\partial} &{}\H^{2}(K,Z)
   }
   \]
in which the map $m_{*}$ is such that
$m(\alpha_{1}, \dots, \alpha_{r}) = \alpha_{1} \dots \alpha_{r}$. 
So, if $\xi \in \H^{1}(K,A)$ is the image of 
$(\xi_{1}, \dots , \xi_{r})$, we have that $\partial\xi$ 
is the class of the product
   \[
   (u_{1}, u_{r+1})_{d_{1}} \otimes_{K} (u_{2}, u_{r+2})_{d_{2}}
      \otimes_{K} \dots \otimes_{K} (u_{r}, u_{2r})_{d_{r}},
   \]
whose index is $d_{1} \dots d_{r}$. Hence 
$\ind(G,Z) \geq d_{1} \dots d_{r} = \sqrt{|A/K|}$, as needed.

Now we prove the claim. Choose a power of $p$, 
call it $d$, that is divisible by the order of $Z$ 
and by the order of each $a_{i}$. 
Consider the group $\Lambda(d)$ defined by the presentation
   \[
   \langle x_1, x_2, y \mid x_1^d = x_2^d = y^d = 1 , \; 
   x_1 x_2 = y x_2 x_1, \; x_1 y = yx_1, \; x_2 y = y x_2 \rangle.
   \]
Call $\rho_{i}\colon \Lambda(d) \arr G_{i}$ the homomorphism obtained 
by sending $x_{1}$ to $a_{i}$, $x_{2}$ to 
$a_{r+i}$, and $y$ to $z^{n/d_{i}} = \omega(a_{i}, a_{r+i})$.

Let $\zeta_{d}$ be a primitive $d$-th root of $1$ 
in $k$ such that $\zeta_{n} = \zeta_{d}^{n/d}$. 
The subgroup $\generate{y}$ in $\Lambda(d)$ is 
cyclic of order $d$; we fix the isomorphism 
$\generate{y} \simeq \mu_{d}$ so that $y$ corresponds 
to $\zeta_{d}$. The restriction of $\rho_{i}$ 
to $\generate{y} \arr Z$ corresponds to 
the homomorphism $\mu_{d} \arr \mu_{n}$ defined 
by $\alpha \arrto \alpha^{d/d_{i}}$. We have a commutative diagram
   \[
   \xymatrix{
   1 \ar[r]
 &\ar[r]\ar[d]_{\substack{\alpha\\\downarrow\\\alpha^{d/d_{i}}}{}} \mu_{d} & 
     {}\Lambda(d)\ar[r]\ar[d]^{\rho_{i}}
     & {}(\ZZ/d\ZZ)^{2} \ar[d]\ar[r] & 1\\
   1 \ar[r] & \mu_{n} \ar[r] & G_{i} \ar[r] & A_{i} \ar[r] & 1\hsmash{.}
   }
   \]
We have $\H^{1}\bigl(K, (\ZZ/d\ZZ)^{2}\bigr) = (K^{*}/{K^{*}}^{d})^{2}$. According to \cite[Example 7.2]{vela}, the image of the element 
$(u_{i}, u_{r+i})\in \H^{1}\bigl(K, (\ZZ/d\ZZ)^{2}\bigr)$ 
is the cyclic algebra $(u_{i}, u_{r+i})_{d}$; 
hence, if $\xi_{i}$ is the image in $\H^{1}(K, A_{i})$ 
of $(u_{i}, u_{r+i})$, the image of $\xi_{i}$ in 
$\H^{2}(K, \mu_{d})$ is the algebra 
$(u_{i}, u_{r+i})_{d}^{\otimes d/d_{i}}$, which is equivalent 
to $(u_{i}, u_{r+i})_{d_{i}}$. This concludes the proof 
of Theorem~\ref{t.p-groups1}(a).
\end{proof}

\begin{proof}[Proof of Theorem~\ref{t.p-groups1}(b)]
By Lemma~\ref{lem.reduction}
we may assume that $\rC(G) = Z$ is cyclic. In this case
Theorem~\ref{t.p-groups1} asserts that $G$ has
a faithful representation of degree~$\sqrt{|A|}$.

Suppose $|\rC(G)| = p^{h}$ and $|A| = p^{2m}$; we want to construct 
a faithful representation of $G$ of dimension~$p^{m}$. 
By \cite[\S 3.1]{tignol-amitsur} $A$ contains a Lagrangian 
subgroup $L$ of order $p^{m}$. 
Denote by $H$ the inverse image of $L$ in $G$; 
then $H$ is an abelian subgroup of $G$ of order $p^{h+m}$. 
Since $\zeta_{p^{e}} \in k$ we can embed $Z$ in $k^{*}$ 
and extend this embedding to a homomorphism 
$\chi\colon H \arr k^{*}$. We claim that 
the representation $\rho\colon G \arr \GL_{p^{m}}$ 
induced by $\chi$ is faithful.

It is enough to show that $\rho(g) \neq \id$ for any $g \in G$ 
of order~$p$, or, equivalently, that 
$\rho\mid_{\generate{g}}$ is non-trivial for any such $g$. 
If $s \in G$ consider the subgroup $H_{s} 
\eqdef s\generate{g}s^{-1} \cap H$ of $H$, which is embedded 
in $\generate{g}$ via the homomorphism 
$x \arrto s^{-1}xs$. By Mackey's formula 
(\cite[\S 7.3]{serre}), $\rho\mid_{\generate{g}}$  
contains all the representations of $\generate{g}$ 
induced by the restrictions $\chi\mid_{H_{s}}$ via the embedding above.

If $g \notin H$ then 
$H_{1} = \generate{g} \cap H = \{1\}$: 
we take $s = 1$, and we see that $\rho\mid_{\generate{g}}$ 
contains a copy of the regular representation 
of $\generate{g}$, which is obviously non-trivial.

Assume $g \in H$. Then $H_{s} = \generate{sgs^{-1}}$ for any $s \in G$; it is enough to prove that $\chi(sgs^{-1}) \neq 1$ for some $s \in G$. If $\chi(g) \neq 1$ then we take $s = 1$. Otherwise $\chi(g) = 1$; in this case $g \notin \rC(G)$, because $\chi\mid_{\rC(G)}\colon  \rC(G) \arr k^{*}$ is injective. Hence the image $\overline{g}$ of $g$ in $A$ is different from $0$, and we can find $s \in G$ such that $\omega(\overline{s}, \overline{g}) \neq 1$. Then we have
\[   \chi(sgs^{-1}) = \chi\bigl(\omega(\overline{s}, \overline{g})g\bigr)
   = \chi\bigl(\omega(\overline{s}, \overline{g})\bigr) \chi(g)
   = \chi\bigl(\omega(\overline{s}, \overline{g})\bigr)
   \neq 1. \]
This concludes the of Theorem~\ref{t.p-groups}(b).
\end{proof}

\section{Essential dimension of $p$-groups II}

We begin with several simple illustrations of Theorem~\ref{t.p-groups}.

\begin{example} \label{ex.extraspecial}
Recall that a $p$-group $G$ is called \emph{extra-special} if its 
center $Z$ is cyclic of order $p$, and the quotient $G/Z$ is elementary 
abelian. The order of an extra special $p$-group $G$ is an odd power 
of $p$; the exponent of $G$ is either $p$ or $p^{2}$; 
cf.~\cite[pp. 145--146]{robinson}.  Note that every non-abelian 
group of order $p^{3}$ is extra-special. For extra-special
$p$-groups Theorem~\ref{t.p-groups} reduces to the following.

{\em Let $G$ be an extra-special $p$-group of order $p^{2m+1}$. 
Assume that that the characteristic of $k$ is different from $p$, 
that $\zeta_{p} \in k$, and $\zeta_{p^{2}} \in k$ if the exponent 
of $G$ is $p^{2}$. Then $\ed G = p^{m}$.}
\end{example}

\begin{example} \label{ex3.p-groups}
Let $p$ be an odd prime and $G = C_{p^r} \ltimes C_{p^s}$ be the natural
semidirect product of cyclic groups of order $p^r$ and $p^s$
(in other words, $C_{p^s}$ is identified with the unique subgroup
of $C_{p^r}^*$ of order $p^s$).  If $s \le r/2$ then
\[ \ed_k(C_{p^r} \ltimes C_{p^s}) = p^s \, , \] 
for any field $k$ containing a primitive $p$th root of unity $\zeta_p$.
\end{example}    

\begin{proof} Here $\rC(G)$ is the (unique) subgroup of $C_{p^r}$ of order
$p^s$. If $s \le r/2$, this subgroup is central. Thus, if $\zeta_{p^r}$,
the equality $\ed_k(G) = p^s$
is an immediate consequence of Theorem~\ref{t.p-groups}. Since
we are only assuming that $\zeta_p \in k$, Theorem~\ref{t.p-groups}
only tells us that $\ed_k(G) \ge p^s$.  
To prove the opposite inequality, we argue as follows.
Let $F$ be the prime subfield of $k$.
By~\cite[Corollary to Proposition 2]{ledet},
$\ed_{F(\zeta_p)}(G) \le p^s$.
Since we are assuming that $F(\zeta_p) \subset k$, 
we conclude that $\ed_k(G) \le p^s$ as well.
\end{proof}

\begin{corollary} \label{cor2.p-groups}
Suppose $k$ is a base field of characteristic $\neq p$.
If $G$ is a non-abelian finite $p$-group then $\ed G \ge p$.
\end{corollary}

\begin{proof}
We argue by contradiction. Assume the contrary and
let $G$ be a non-abelian $p$-group $G$ of smallest possible order 
such that $\ed G < p$. Since $G$ has a non-trivial center,  
there exists a cyclic central subgroup $Z \subset G$.
The short exact sequence 
\begin{equation} \label{e.exact-sequence}
1  \arr Z  \arr G  \arr G/Z  \arr 1
\end{equation}
give rise to the exact sequence of pointed sets
   \[
   \H^1(K, G) \arr \H^1(K, G/Z) \stackrel{\partial_K}\larr \H^2(K, Z)
   \]
for any field extension $K$ of our base field $k$.
We will now consider two cases.

{\bf Case 1.} Suppose the map $\H^1(K, G) \arr \H^1(K,
G/Z)$ is not surjective for some $K/k$. Then $\partial_K$ is
non-trivial, and Theorem~\ref{thm2b} tells us that $\ed G \ge p$, a
contradiction.

{\bf Case 2.} Suppose the map $\H^1(K, G)  \arr \H^1(K, G/Z)$ is 
surjective for every $K/k$. Then
the morphism $\cB G \arr \cB(G/Z)$ is isotropic, and 
Proposition~\ref{p.isotropic} implies that
$p > \ed G \ge \ed(G/Z)$. By the minimality of $G$, the group 
$G/Z$ has to be abelian. Consequently, $[G, G] \subset Z$ is cyclic 
and central in $G$. Since $G$ is non-abelian, $|G/\rC(G)| \ge p^2$. 
Theorem~\ref{t.p-groups} now tells us that 
\[ \ed(G) = \sqrt{|G|/|\rC(G)|} + \rank \, \rC(G) - 1 \ge p \, , \] 
a contradiction.
\end{proof}

We will conclude this section by answering the following question
of Jensen, Ledet and Yui~\cite[p. 204]{jly}.

\begin{question}
\label{q:ljy} Let $G$ be a finite group and $N$ be a normal subgroup.
Is it true that $\ed G \ge \ed(G/N)$?
\end{question}

The inequality $\ed(G) \ge \ed(G/N)$ is known to hold in many 
cases (cf., e.g., Lemma~\ref{lem5.3}). We will now show that 
it does not hold in general, even if $H$ is assumed to be central.

\begin{corollary} \label{cor.jly}
For every real number $\lambda > 0$ there exists a finite $p$-group $G$, 
with a central subgroup $H \subset G$ such that $\ed(G/H) > \lambda \ed G$.
\end{corollary}

\begin{proof} Let $\Gamma$ be a non-abelian group of order $p^3$.
The center of $\Gamma$ has order $p$; denote it by $C$. The
center of $\Gamma^n = \Gamma \times \dots \times \Gamma$ 
($n$ times) is then $C^n$. Let $H_n$ be the subgroup of $C^n$ 
consisting of $n$-tuples $(c_1, \dots, c_n)$ such that 
$c_1 \dots c_n = 1$.  Clearly 
   \[
   \ed\Gamma^n \le n \cdot \ed\Gamma = np \, ;
   \]
see Example~\ref{ex.extraspecial}.

On the other hand, $G_n/H_{n}$, is easily seen to be extra-special 
of order $p^{2n+1}$, so $\ed(G_{n}/H_{n}) = p^{n}$, again by 
Example~\ref{ex.extraspecial}. Hence by taking $n$ sufficiently 
large we finish the proof.
\end{proof}

\section{Spinor groups}
\label{s.spinor}

In this section we will prove Theorem~\ref{thm.spin} stated in the
introduction. For an introduction to the structure of this group
and the theory of Clifford algebras and spin modules, we refer the reader
to \cite[\S 20.2]{fulton-harris} and \cite{chevalley}.

Recall that by $\Spin_{n}$ we mean the totally split form of the spin
group in dimension~$n$ over a field of characteristic not equal to
$2$.  Let us be explicitly about this.

Write $\langle a_1,\ldots, a_n \rangle$ for the rank $n$-quadratic form $q$
given by $q(x_1,\ldots, x_n)=\sum_{i=1}^n a_i x_i^2.$   Set 
$h$ to be the standard hyperbolic quadratic form given by
$h(x,y)=xy$.  (Thus $h\cong \langle 1, -1 \rangle$).  For each
$n\ge 0$ define
\begin{equation}
h_n =
\begin{cases}
  h_n^{\oplus n/2}, & \text{if $n$ is even,}\\
  h_n^{\oplus (n-1/2)} \oplus \langle 1 \rangle, & \text{if $n$ is odd.}
\end{cases}
\end{equation}
Set $\Orth_n\eqdef\Orth(h_n), \SO_n\eqdef\SO(h_n)$, and 
$\Spin_n=\Spin(h_n)$.  These are all totally split groups.

Now, one of the hypotheses of Theorem~\ref{thm.spin} is that
$\zeta_{4} \in k$.
Therefore we can write $\Spin_n$ 
as $\Spin(q)$, where
   \[
   q(x_{1}, \dots, x_{n}) = -(x_{1}^{2} + \dots + x_{n}^{2}).
   \]

Consider the subgroup $\Gamma_{n} \subseteq \SO_{n}$ consisting 
of diagonal matrices, which is isomorphic to 
$\mu_{2}^{n-1}$. Call $G_{n}$ the inverse image of 
$\Gamma_{n}$ in $\Spin_{n}$. It is a constant group scheme over $k$.  
Denote by $\mu_{2}$ the kernel 
of the homomorphism $\Spin_{n} \arr \SO_{n}$.

\begin{lemma}
Every $\Spin_{n}$-torsor over an extension $K$ of $k$ admits a reduction of structure group to $G_{n}$.
\end{lemma}

\begin{proof}
Let $P \arr \Spec K$ be a $\Spin_{n}$-torsor: we are claiming that 
the $K$-scheme $P/G_{n}$ has a rational point. We have 
$P/G_{n} = (P/\mu_{2})/\Gamma_{n}$. However 
$P/\mu_{2} \arr \Spec K$ is the $\SO_{n}$ torsor associated 
with $P \arr \Spec K$, and every $\SO_{n}$-torsor has 
a reduction of structure group to $\Gamma_{n}$.
\end{proof}

This means that the natural morphism $\cB G_{n} \arr \cB\Spin_{n}$ 
is isotropic; so from Propositions \ref{p.RelDim-quotients} 
and \ref{p.isotropic} we get the bounds
   \begin{equation}\label{eq:bound-spin}
   \ed G_{n} - \dim\Spin_{n} \leq \ed\Spin_{n} \leq \ed G_{n}.
   \end{equation}
Of course $\dim \Spin_{n} = n(n-1)/2$; we need to compute 
$\ed G_{n}$. The structure of $G_{n}$ is well understood 
(besides the references cited above, 
it is also very clearly described in \cite{wood}). 
The group scheme $\Spin_{n}$ is a subgroup scheme of the group 
scheme of units in the Clifford algebra $A_{n}$ of the quadratic 
form $-(x_{1}^{2} + \dots + x_{n}^{2})$. The algebra $A_{n}$ 
is generated by elements $e_{1}$, \dots,~$e_{n}$, 
with relations $e_{i}^{2} = -1$ and 
$e_{i}e_{j} + e_{j}e_{i} = 0$ for all $i \neq j$. 
The element $e_{i}$ is in $\Pin_{n}$, and image 
of $e_{i}$ in $\rO_{n}$ is the diagonal matrix with 
$-1$ as the $i$-th diagonal entry, and $1$ as all the other 
diagonal entries. The kernel of the homomorphism 
$\Pin_{n} \arr \rO_{n}$ is $\{\pm 1\}$.

For any $I \subseteq \{1, \dots, n\}$ 
write $I = \{i_{1}, \dots, i_{r}\}$ with 
$i_{1} < i_{2} < \dots < i_{r}$ and set 
$e_{I} \eqdef e_{i_{1}} \dots e_{i_{r}}$. 
The group $G_{n}$ consists of the elements of 
$A_{n}$ of the form $\pm e_{I}$, 
where $I \subseteq \{1, \dots, n\}$ has an even number 
of elements. The element $-1$ is central, and 
the commutator $[e_I, e_J]$ is given by
   \[
   [e_{I}, e_{J}] = (-1)^{|I \cap J|}
   \]
It is clear from this description that $G_n$ is a $2$-group,
of order $2^n$, the commutator $[G_n, G_n] = \{ \pm 1 \}$ is cyclic,
and the center $\rC(G)$ is given by 
\[ \rC(G_n) = \begin{cases} 
\text{$\{ \pm 1 \} \simeq \ZZ/2 \ZZ$, if $n$ is odd,} \\
\text{$\{ \pm 1, \pm e_{\{ 1, \dots, n \}} \} \simeq \ZZ/4 \ZZ$, 
if $n \equiv 2$ (mod $4$),} \\
\text{$\{ \pm 1, \pm e_{ \{ 1, \dots, n \} } \} \simeq \ZZ/2 \ZZ \times \ZZ/2 \ZZ$, 
if $n$ is divisible by $4$.} \end{cases} \]
Theorem~\ref{t.p-groups} now tells us that
\[ \ed(G_n) = \begin{cases} \text{$2^{(n-1)/2}$, if $n$ is odd,} \\ 
 \text{$2^{(n-2)/2}$, if $n \equiv 2$ (mod $4$),} \\
 \text{$2^{(n-2)/2} + 1$, if $n$ is divisible by $4$.} \end{cases} \]
Substituting this into~\eqref{eq:bound-spin}, we obtain the bounds of
Theorem~\ref{thm.spin}.
\qed

\begin{remark} \label{rem.pin} The same argument, with $G_n$ replaced by
the inverse image of the diagonal subgroup of $O_n$ in $\Pin_n$, yields 
the following bounds on the essential dimensions of $\Pin$ groups
(over a field $k$ satisfying the assumptions of Theorem~\ref{thm.spin}): 
\begin{eqnarray*} 
2^{\lfloor n/2 \rfloor} - \frac{n(n-1)}{2} \le
   \ed \Pin_{n} \le  2^{\lfloor n/2 \rfloor}, & 
\text{if $n \not\equiv 1 \pmod 4$,} \\
2^{\lfloor n/2 \rfloor} - \frac{n(n-1)}{2} + 1 \le
   \ed \Pin_{n} \le  2^{\lfloor n/2 \rfloor} + 1, &
\text{if $n \equiv 1 \pmod 4$.} \end{eqnarray*}
\end{remark}

\begin{remark}
When $n \leq 14$ the lower bound of Theorem~\ref{thm.spin} is 
negative and the upper bound is much larger than the true value 
of $\ed\Spin_{n}$. For $n = 15$ and $16$ our inequalities yield
   \[
   23 \leq \ed \Spin_{15} \leq 128
   \]
and
   \[
   9 \leq \ed \Spin_{16} \leq 129.
   \]
When $n = 16$ our lower bound coincides with the lower 
bound~\eqref{e.old-spin} of Reichstein--Youssin and 
Chernousov--Serre, while for $n = 15$ it is substantially larger. 
When $n \geq 17$ the exponential part of the lower bound takes 
over, the growth becomes fast and the gap between the lower 
bound and the upper bound proportionally small. For values 
of $n$ close to $15$ our estimates are quite imprecise; 
it would be interesting to improve them.
\end{remark}

\begin{remark}
By Proposition~\ref{p.extensions}, the lower bounds in the theorem hold for 
over any field of characteristic different from $2$ (and for any form 
of the Spin group).

On the other hand, if we do not assume that $\zeta_{4} \in k$,
we get the slightly weaker upper bound
   \[
   \ed\Spin_{n} \leq 2^{\lfloor(n-1)/2\rfloor} + n - 1
   \]
for the totally split form of the spin group in dimension~$n$.
To prove this inequality, we observe that
a generically free representation of $\Spin_{n}$ 
can be constructed 
by taking a spin, or half-spin, representation $V$ of 
$\Spin_{n}$ of dimension $2^{\lfloor(n-1)/2\rfloor}$, 
and adding a generically free representation $W$ of $\SO_{n}$. 
Since the essential dimension of $\SO_{n-1}$ 
is $n-1$ over any field of characteristic different from~$2$,
there is an $\SO_n$-compression $f \colon W \darr X$, where $\dim(X) = 
\dim(\SO_n) + n-1$. Now
$\id \times f \colon V \times W \darr V \times X$ is 
a $\Spin_n$-compression of $V \times W$.  Consequently, 
\begin{align*}
 \ed \Spin_n & \ge \dim(V \times X) - \dim\Spin_n \\
 & =  2^{\lfloor(n-1)/2\rfloor} + \dim\SO_n + n - 1 - \dim\Spin_n \\
 & = 2^{\lfloor(n-1)/2\rfloor} + n - 1 \, , \end{align*}
as claimed.
\end{remark}

\begin{remark} \label{rem.clifford}
It is natural to ask whether the inequality
\[ \ed \Spin_n  \ge 2^{\lfloor (n -1)/2 \rfloor} - \frac{n(n-1)}{2} \]
can be proved by a direct application of Theorem~\ref{thm.spin}
to the exact sequence 
\begin{equation}
\label{e.SpinExact}
 1 \arr \mu_2 \arr \Spin_n \arr \SO_n \arr 1 \,  
\end{equation}
without considering the finite subgroup $G_n$ of $\Spin_n$.
The answer is ``yes.'' Note, however, that we see no way of obtaining
the slightly stronger lower-bound which we obtain when $n$ is
divisible by $4$.

Indeed, consider the associated coboundary map
 \[ \xymatrix{ \H^1(K, \SO_{m}) \ar@{->}[r]^{\partial_K} & \H^2(K, \mu_2).} \]
A class in $\H^1(K, \SO_m)$ is represented by a $m$-dimensional
quadratic form $q$ of discriminant $1$ defined over $K$. The
class of $\partial_K(q) \in H^2(K, \mu_2)$ is then
the Hasse-Witt invariant of $q$; following Lam~\cite{lam},
we will denote it by $c(q)$. (Note that since we are assuming that
$-1$ is a square in $k$, the Hasse invariant and the Witt
invariant coincide; see \cite[Proposition V.3.20]{lam}.)
Our goal is thus to show that for every $n \ge 1$ there exists 
a quadratic form $q_n$ of dimension $n$ and discriminant $1$ such 
that $c(q_n)$ has index $2^{\lfloor (n -1)/2 \rfloor}$.

If $n$ is even this is proved in~\cite[Lemma 5]{m2}.
(Note that in this case $c(q) \in H^2(K, \mu_2)$ is the class of               
the Clifford algebra of $q$.) If $n = 2m + 1$ is odd, set 
$K = k(a_1, b_1, \dots, a_m, b_m)$, where 
$a_1, b_1, \dots, a_m, b_m$ are independent variables,
and define $q_n$ recursively by                      
\[ \text{$q_3 = \langle a_1, b_1, a_1b_1 \rangle$ and 
$q_{n + 2} = \langle a_n b_n \rangle \otimes q_n \oplus
\langle a_n, b_n \rangle$.} \]   
A direct computation using basic properties of the Hasse-Witt invariant
(see, e.g., \cite[Section V.3]{lam}) shows that 
$c(q_{2m + 1})$ is the class of the product 
$(a_1, b_1)_2 \otimes_K \dots \otimes_K (a_m, b_m)_2$ of quaternion algebras.
This class has index $2^m$, as claimed.
\qed

Any lower bound on $\ed \Spin_n$ obtained in this way
(i.e., directly from Theorem~\ref{thm2b}) will necessarily 
be of the form $\ed \Spin_n  \ge 2^m - \frac{n(n-1)}{2}$ for some 
integer $m \ge 0$.  Thus, while this approach recovers the lower 
bound of Theorem~\ref{thm.spin} if $n$ is not divisible by $4$,
it cannot be used to do so if $n$ is divisible by $4$. 
\end{remark}

To conclude this section, we will now prove similar bounds on the essential 
dimensions on half-spin groups. We begin with the following
simple corollary of~\cite[Theorem 1.1]{cgr}, 
which appears to have been previously overlooked.

\begin{lemma} \label{lem.cgr}
Let $G$ be a closed (but not necessarily connected)
subgroup of $\GL_{n}$ defined over a field $k$. Assume 
that one of the following conditions holds.

\begin{enumerate}

\item $k$ is algebraically closed of characteristic~$0$.

\item $k$ has characteristic~$0$ and $G$ is connected.

\item $G$ is connected and reductive.

\end{enumerate}

Then $\ed G \leq n$.
\end{lemma}

\begin{proof}
According to \cite[Theorem~1.1]{cgr}, there exists a finite subgroup 
scheme $S \subseteq G$ such that every $G$-torsor over $\spec K$, 
where $K$ is an extension of $G$, has a reduction of structure groups 
to $S$. Hence the morphism $\cB S \arr \cB G$ is isotropic, 
so $\ed G \leq \ed S$. But the restriction of the representation 
$G \subseteq \GL_{n}$ is generically free, hence $\ed S \leq n$.
\end{proof}

\begin{example} \label{ex.cgr1}
Suppose $G$ satisfies one of the conditions (a), (b) 
or (c) of Lemma~\ref{lem.cgr} and the centralizer $C_G(G^0)$
of the connected component of $G$ is trivial. Then 
the adjoint representation of $G$ is faithful and
Lemma~\ref{lem.cgr} tells us that $\ed(G) \le \dim(G)$.
In particular, this inequality is valid for every connected semisimple 
adjoint group $G$. (In the case of simple adjoint groups, a stronger 
bound is given by \cite[Theorem 1.3]{lemire}.)
\end{example}

We are now ready to proceed with our bounds on the essential 
dimension of half-spin groups. Recall that 
the \emph{half-spin group} $\HSpin_{n}$ is defined, for every $n$ 
divisible by $4$, as $\Spin_{n}/\generate{\eta}$, where $\eta$ 
is an element of the center of $\Spin_{n}$ different from $-1$. 
(There are two such elements, but
the resulting quotients are isomorphic.)

\begin{theorem}
Suppose $k$ is a field of characteristic $\neq 2$ and 
$\zeta_{4} \in k$. Let $n$ be a positive integer divisible by $4$. Then
   \[
   2^{(n-2)/2} - \frac{n(n-1)}{2} \leq \ed\HSpin_{n} \leq 2^{(n-2)/2}
   \]
\end{theorem}

\begin{proof}
The group $\HSpin_{n}$ contains $G_{n}/\generate{\eta} \simeq G_{n-1}$, 
which is an extra-special group of order $2^{n-1}$. By 
Example~\ref{ex.extraspecial}
$\ed(G_{n}/\generate{\eta}) = 2^{(n-2)/2}$ and thus
   \begin{align*}
   \ed \HSpin_n &\ge \ed(G_{n}/\generate{\eta}) - \dim \HSpin_n \\
   &= 2^{(n-2)/2} - \frac{n(n-1)}{2} \, ,
   \end{align*}
as in the proof of Theorem~\ref{thm.spin}.

For the upper bound notice that one of the two half-spin 
representations of $\Spin_{n}$ descends to $\HSpin_{n}$, 
and is a faithful representation of $\HSpin_{n}$ of dimension
$2^{(n-2)/2}$. The upper bound now follows from Lemma~\ref{lem.cgr}
\end{proof}

\section{Essential dimension of cyclic $p$-groups}
\label{s.cyclic}

In this section, we are going to prove the following theorem
due to M.~Florence.
In the sequel $p$ will denote a prime, different from the characteristic 
of our base field $k$, and 
$\zeta_d$ will denote a primitive $d$th root of unity in $\overline{k}$. Recall that we have set $C_{n} \eqdef \ZZ/n\ZZ$.

\begin{theorem}[M.~Florence {\cite{florence}}]
\label{t.cyclic}
Let $p$ be a prime, $k$ a field of characteristic $\neq p$.
Suppose $\zeta_{p^n} \in k$ but $\zeta_{p^{n+1}}  \notin k$ for some 
integer $n \ge 1$.  Moreover, if $p = 2$ and $n = 1$, assume also
that $k(\zeta_4) \neq k(\zeta_8)$. Then
$$ \ed C_{p^m} =
\begin{cases}
p^{m-n} &  \text{if }n<m,\\
1   &  \text{if }n\geq m.
\end{cases}
$$
for any integer $m\geq 1$, 
\end{theorem}

This theorem was independently obtained by us in the case where 
$n \ge \lfloor(m + 1)/2\rfloor$; in particular, for $m = 2$. 
However, our proof of the stronger result given by 
Theorem~\ref{t.cyclic}, will rely on 
an idea of M.~Florence, in combination with  
the lower bound of Theorem~\ref{thm2b}.

\begin{proof} If $m\leq n$, then $C_{p^m}=\mu_{p^m}$.  Therefore $\ed
  C_{p^m} =1$ by~\cite[Example 2.3]{bf1}.  We can therefore restrict
  our attention to the case $n<m$.

  We first show that $\ed C_{p^m} \leq p^{m-n}$.  To do this, pick a
  faithful character $\chi\colon  C_{p^n} \arr\GG_m$ defined over $K$ and set
  $V\eqdef \Ind_{C_{p^n}}^{C_{p^m}} \chi$.  A simple calculation shows
  that $V$ is faithful, thus, $V$ is generically free since $C_{p^m}$
  is finite.  By Theorem~\ref{t.GenFree}, it follows that $\ed
  C_{p^m} \leq \dim V=p^{m-n}$.

  To show that $\ed C_{p^m}\geq p^{m-n}$, view the representation $V$
as a homomorphism $\rho\colon C_{p^m} \arr \GL(V)$ of algebraic groups.  Let
$\pi:\GL(V) \arr \PGL(V)$ denote the obvious projection and note that
the kernel of $\pi\circ\rho$ is exactly $C_{p^n}$. It follows that we
have a commutative diagram
$$
\xymatrix{ 0\ar[r] & C_{p^n}\ar[r]\ar[d] & C_{p^m}\ar[r]^{\rho}\ar[d]
& C_{p^{m-n}}\ar[r]\ar[d]^{\iota} & 1\\ 1\ar[r] & \GG_m\ar[r]
&\GL(V)\ar[r]^{\pi} &\PGL(V)\ar[r] & 1 \\ }
$$
where the rows are exact and the columns are injective.

Let $K/k$ be a field extension and let $t\in \H^1(K,C_{p^m})$ be a
torsor.  Let $\iota_*\colon \H^1(K,C_{p^{m-n}}) \arr \H^1\bigl(K, \PGL(V)\bigr)$
denote the map induced by $\iota$.  Then, from the commutativity of
the above diagram (and the injectivity of the columns), it follows
that $\ind_K(t)$ is the index of of the CSA $\iota_*(t)$.

We claim that there is a $t$ such that $\iota_*(t)$ is a division
algebra. From this it will clearly follow that $\ind_K(t)=\dim
V=p^{m-n}$.

In fact, this $t\in \H^1(K,C_{p^m})$ is simply the ``generic'' one.
(This is the part of the argument that we learned from 
Florence's preprint~\cite{florence}.)
Namely, let $L=K(x_1,\ldots, x_{p^{m-n}})$ denote the field obtained
by adjoining $p^{m-n}$ independent variables to $K$, and let
$C_{p^{m-n}}$ act on $L$ by permuting the variables in the obvious way
($k\cdot x_i=x_{i+k}\pmod{p^{m-n}}$).  Let $F=L^{C_{p^{m-n}}}$.  Then $L/K$
  defines a $C_{p^{m-n}}$-torsor $t$ over $K$.  

In the case where $k = \QQ(\zeta_{p^n})$
the torsor $\iota_*(t)$ is the ``generic'' algebra $R_{p^n,
p^m, p^m}$ of~\cite[\S 7.3]{rowen2}.  By a theorem of Brauer (see 
Theorem~\cite[Theorem~7.3.8]{rowen2}) it is a division algebra.
A similar argument (due to M. Florence) shows that the same is true
if $\QQ(\zeta_{p^n})$ is replaced by
our field $k$ (satisfying the assumptions of Theorem~\ref{t.cyclic}).
\end{proof}

\begin{remark} \label{rem36} Suppose $\zeta_6 \in k$.
Note that in this case $[k(\zeta_{36}:k]$ always divides $6$. 
Assume $[k(\zeta_{36}): k] = 6$. (This occurs, for example, 
if $k = \QQ(\zeta_6)$.) We claim that $\ed_k C_{36} \ge 4$.

Indeed, by Remark~\ref{rem.ctkm} it suffices to show 
that $\ind(C_{36}, C_6) = 6$.  Let $K/k$ be a field extension 
and consider the boundary map 
   \[
   \partial_K \colon K^*/(K^*)^6  = \H^{1}(K, C_6) \arr \H^{2}(K, C_6)
   \]
induced by the exact sequence
$1 \to C_6 \to C_{36} \to C_6 \to 1$.
By \cite[Theorem 7.1]{vela}, $\partial_K$ sends $(a)$ to 
the class of the cyclic algebra $(a, \zeta_6)_6$. The index 
of this cyclic algebra clearly divides $6$. Taking $K = k(a)$, 
where $a$ is an independent variable over $k$, and applying
Wedderburn's criterion (cf. e.g.,
\cite[Corollary 15.1d]{pierce}), 
we conclude that in this case the cyclic algebra $(a, \zeta_6)_6$ has
index $6$. Thus $\ind(C_{36}, C_6) = 6$, as claimed.

A similar argument shows that if Conjecture~\ref{conj:ed-product}
is valid for $n = p_{1}^{a_{1}} \dots p_{r}^{a_{r}}$ then
$ \ed_k(C_n) \ge p_{1}^{a_{1}} + \dots + p_{r}^{a_{r}} - r + 1$.
\end{remark}

\smallskip
Let $D_n$ be the dihedral group of order $2n$.
Ledet~\cite[Section 3]{ledet} conjectured that if $n$ is odd
then $\ed C_n = \ed D_n$ over any field $k$ of characteristic zero. 
We will now prove
this conjecture in the case where $n = p^r$ is a prime power
and $k$ contains a primitive $p$th root of unity.

\begin{corollary} \label{cor.ledet}
Let $p$ be an odd prime and
$k$ be a field containing a primitive $p$th root of unity.
Then $\ed_k D_{p^m} = \ed_k C_{p^m}$.
\end{corollary}

\begin{proof}  If $\zeta_{p^m} \in k$ then we know that
$\ed_k C_{p^m} = \ed_k D_{p^m} = 1$;  see the proof
of~\cite[Theorem 6.2]{bur}.  Thus we may assume $\zeta_{p^m} \not \in k$.

Let $s$ be the largest integer $n$ such that
$\zeta_{p^n} \in k$. By our assumption $1 \le s \le m-1$.
By Theorem~\ref{t.cyclic}
   \[
   \ed_k C_{p^m} = p^{m - n} \, .
   \]
Since $C_{p^m} \subset D_{p^m}$, we clearly
have $\ed_k D_{p^m} \ge \ed_k C_{p^m}$. Thus we only need to show that
\[ \ed_k D_{p^m} \le  p^{m-n} \, . \]
To prove this inequality, note that
$D_{p^m} \simeq C_{p^m} \rtimes C_2$
has a subgroup isomorphic to $D_{p^n} = C_{p^n} \rtimes C_2$
of index $p^{m-n}$. Since $k$ contains $\zeta_{p^n}$,
$D_{p^n}$ has essential dimension $1$ over $k$.
Thus, by \cite[Section 3]{ledet},
\[ \ed_k D_{p^m} = \ed_k D_{p^n} \le [D_{p^m} : D_{p^n}]
 = 1 \cdot p^{m - n} = p^{m - n}. \]
This completes the proof of Corollary~\ref{cor.ledet}.
\end{proof}

\section{Pfister numbers}
\label{s.arason}

Let $k$ be a field of characteristic not equal to $2$ and write
$\rW(k)$ for the Witt ring of $k$; see~\cite[Chapter 2]{lam}.  Let
$I=I(k)$ denote the ideal of all even dimensional forms in the Witt
ring.  Then, for any integer $a>0$, $I^a$ is generated as an abelian
group by the $a$-fold Pfister forms~\cite[Proposition 1.2]{lam}.

Let $q$ be a quadratic form of rank $n>0$ whose class $[q]$ in $\rW(k)$
lies in $I^a$ for $a>0$.  Define the \emph{$a$-Pfister number} of $q$ to be
the minimum number $r$ appearing in a representation 
$$
q=\sum_{i=1}^r \pm p_i
$$
with the $p_i$ being $a$-fold Pfister forms.  The
\emph{$(a,n)$-Pfister number} $\Pf_k(a, n)$ is the supremum of the
$a$-Pfister number of $q$ taken over all field extensions $K/k$ and 
all $n$-dimensional forms $q$ such that $[q] \in I^a(K)$.   

We have the following easy (and probably well-known) result.

\begin{proposition}
\label{p.PfisterEasy}  
Let $k$ be a field of characteristic not equal to
  $2$ and let $n$ be a positive even integer. 
  \begin{enumeratea}
  \item 
  $\Pf_k(1,n) \leq n$.
  \item 
  $\Pf_k(2,n) \leq n - 2$.
  \end{enumeratea}
\end{proposition}
  
\begin{proof}
  (a) If $n$ is even $\langle a_1,\ldots, a_n
  \rangle = \sum_{i=1}^n (-1)^i \da{-a_i}$.  

  (b) Let $q=\langle a_1,\ldots, a_n\rangle$ be an $n$-dimensional
quadratic form over $K$. Recall that 
  $q \in I^2(K)$ iff $n$ is even and $d_{\pm}(q)=1$, modulo 
$(K^*)^2$~\cite[Corollary II.2.2]{lam}. Here $d_{\pm}(q)$ 
is the \emph{signed determinant}
  given by $(-1)^{n(n-1)/2}d(q)$ where $d(q)=\prod_{i=1}^n a_n$ is the
  determinant~\cite[p.38]{lam}.  

To explain how to write $q$ as a sum of $n-2$ Pfister forms,
we will temporarily assume that $\zeta_4 \in K$. In this case 
we may assume that $a_1 \dots a_n = 1$. Since $\langle a, a \rangle$
is hyperbolic for every $a \in K^*$, we see that
$q = \langle a_1, \dots, a_n \rangle$ is Witt equivalent to
\[ \ll a_2, a_1 \gg \oplus 
\ll a_3, a_1a_2 \gg \oplus \cdots \oplus \ll a_{n-1}, a_1 \dots a_{n-2} \gg \, . \]
By inserting appropriate powers of $-1$, we can modify this formula so that
it remains valid even if we do not assume that $\zeta_4 \in K$, as follows:
\[ q = \langle a_1, \dots, a_n \rangle \simeq
\sum_{i=2}^n (-1)^i \da{(-1)^{i+1} a_i, (-1)^{i(i-1)/2 +1} a_1 \dots a_{i-1}} 
\qedhere \]
\end{proof}

We do not have an explicit upper bound on $\Pf_k(3, n)$; however, we do know that $\Pf_k(3, n)$
is finite for any $k$ and any $n$.

To explain this, let us recall that $I^3(K)$ is the set of all classes
$[q]\in\rW(K)$ such that $q$ has even dimension, trivial signed determinant
and trivial Hasse-Witt invariant~\cite{invol}.  

Let $n$ be a positive integer. Let $q$ be a
non-degenerate $n$-dimensional quadratic form over $K$ whose whose signed
determinant is $1$. The class of $q$ in $\H^{1}(K, \Orth_{n})$ lies in
$\H^{1}(K, \SO_{n})$. We say that $q$ \emph{admits a spin structure}
if its class is in the image of $\H^{1}(K, \Spin_{n})$ into $\H^{1}(K,
\SO_{n})$. As pointed out in Remark~\ref{rem.clifford}, the
obstruction to admitting a spin structure is the Hasse-Witt invariant
$c(q)$.  Thus, the forms in $I^3$ are exactly the even dimensional
forms admitting a spin structure. The following result was suggested to us by Merkurjev and Totaro.

\begin{proposition}
\label{p.PfisterEasy2} Let $k$ be a field of characteristic different from~$2$.
  Then $\Pf_k(3,n)$ is finite.
\end{proposition}
\begin{proof}[Sketch of proof]
  Let $E$ be a versal torsor for $\Spin_n$ over a field extension
  $L/k$; cf. \cite[Section I.V]{gms}.
  Let $q_L$ be the quadratic form over $L$ corresponding to
  $E$ under the map $\H^1(L,\Spin_{n})\to \H^1(L,\Orth_{n})$.  The
  $(3,n)$-Pfister number of $q_L$ is then an upper bound for the
  $(3,n)$-Pfister number of any form over any field extension $K/k$.
\end{proof}

\begin{remark}
  For $a > 3$ the finiteness of $\Pf_k(a,n)$ is an open problem.
\end{remark}

The main theorem in this section is a lower bound for 
$\Pf_k(3,n)$ stated as Theorem~\ref{t.pfister} in the Introduction.
We restate it here for the reader's convenience.

\begin{theorem} \label{thm.pfister}
  Let $k$ be a field of characteristic different from $2$
  and let $n$ be an even positive integer. Then 
%
\[ \Pf_k(3,n) \ge \frac{2^{(n+4)/4} -n-2}{7}\, . \]
%
\end{theorem}

For each extension $K$ of $k$, denote by $\rT_{n}(K)$ the image of
$\H^{1}(K, \Spin_{n})$ into $\H^{1}(K, \SO_{n})$; we get a functor
$\rT_{n}\colon \Fields_{k} \arr \Sets$. The essential dimension of
this functor is closely related to the essential dimension of
$\Spin_{n}$.

\begin{lemma} \label{lem.T-Spin} 
$\ed\Spin_{n} - 1 \leq \ed\rT_{n} \leq \ed\Spin_{n}$.
\end{lemma}
\begin{proof}
  In the language off~\cite[Definition 1.12]{bf1}, we 
have a fibration of functors
\[
\H^1(-, \mu_{2})\leadsto \H^1(-,\Spin_n)\arr \rT_{n}(K).
\]
The first inequality then follows from~\cite[Proposition 1.13]{bf1}
and the second follows from Proposition~\ref{p.isotropic}.
\end{proof}

\begin{lemma}\label{lem:sum-spin}
  Let $q$ and $q'$ be non-degenerate quadratic forms over $K$. Suppose
  that $q$ admits a spin structure. Then $q \oplus q'$ admits a spin
  structure if and only if $q'$ admits a spin structure.
\end{lemma}
\begin{proof}
Immediate from the identity $c(q\oplus q')=c(q)+c(q')$~\cite[V.3.15]{lam}.
\end{proof}

Let $h_{K}$ be the standard $2$-dimensional hyperbolic form
$h_{K}(x,y) = xy$ over an extension $K$ of $k$ discussed at the beginning
of \S 13. For each
$n$-dimensional quadratic form $q$ admitting a spin structure over
$K$, denote $\ed_{n}(q)$ the essential dimension of the class of $q$
in $\rT_{n}(K)$.

\begin{lemma}\label{lem:decrease-ed}
Let $q$ be an dimensional quadratic form over $K$ admitting a spin structure, and let $s$ be a positive integer. Then we have
   \[
   \ed_{n+2s}(h_{K}^{\oplus s} \oplus q)\geq \ed_{n}(q)
      - \frac{s(s + 2n -1)}{2}.
   \]
\end{lemma}

\begin{proof}
  Set $m \eqdef \ed_{n+2s}(h_{K}^{\oplus s} \oplus q)$; let $F$ be a
  field of definition of $h^{\oplus s}_{K} \oplus q$ of transcendence
  degree~$m$, and let $\widetilde{q}$ be a quadratic form
  with a spin structure on $F$ whose base change to $K$ is isomorphic
  to $h_{K}^{\oplus s} \oplus q$. Let $X$ be the Grassmannian of
  $s$-dimensional subspaces of $F^{n+2s}$ which are totally isotropic
  with respect to $\widetilde{q}$; the dimension of $X$ is precisely
  $s(s + 2n -1)/2$.

  The variety $X$ has a rational point over $K$; hence there exists an
  intermediate extension $F \subseteq E \subseteq K$ such that
  $\trdeg_{F}E \leq s(s + 2n -1)/2$, with the property that
  $\widetilde{q}_{E}$ has a totally isotropic subspace of
  dimension~$s$. Then $\widetilde{q}_{E}$ splits as $h_{E}^{s} \oplus
  q'$. By Witt's cancellation Theorem, the base change of $q'$ to $K$
  is isomorphic to $q$; hence $\ed_{n}(q) \leq m + s(s + 2n -1)/2$, as
  claimed.
\end{proof}

\begin{proof}[Proof of Theorem~\ref{thm.pfister}]
If $n \leq 10$ then the statement 
is vacuous, because then $2^{(n+4)/4}-n-2 \leq 0$, 
so we assume that $n \geq 12$. We may also assume 
without loss of generality that
$\zeta_{4} \in k$. In this case $\rW(K)$ is a
$\ZZ/2$-vector space; it follows that the $3$-Pfister number 
of a form $q$ is the smallest $r$ appearing in an expression
\[
q=\sum_{i=1}^r \da{a_i,b_i,c_i}.
\]  
in $\rW(K)$.

We will take an $n$-dimensional form $q$ with 
a spin structure such that $\ed_{n}(q) = \ed\rT_{n}$. Suppose 
that $q$ is equivalent in the Witt ring to a form of the type
$\sum_{1 = 1}^{r}\da{a_{i}, b_{i}, c_{i}}$.

Let us write a Pfister form $\da{a,b,c}$ as
   \[
   \da{a,b,c} = \generate{1} \oplus \da{a,b,c}_{0},
   \]
where
   \[
   \da{a,b,c}_{0} \eqdef
   \langle a_i, b_i, c_i, a_ib_i, a_i c_i, b_i c_i, a_i b_i c_i \rangle.
   \]
Set
   \[
   \phi \eqdef \sum_{1 = 1}^{r}\da{a_{i}, b_{i}, c_{i}}_{0}
   \]
if $r$ is even, and
   \[
   \phi \eqdef \generate{1} \oplus \sum_{1 = 1}^{r}\da{a_{i}, b_{i}, c_{i}}_{0}
   \]
if $r$ is odd. Then $q$ is equivalent to 
$\phi$ in the Witt ring, and $\phi$ has 
a spin structure. The dimension of $\phi$ is $7r$ 
or $7r+1$, according to the parity of $r$. 

We claim that $n < 7r$. If not, then the dimension of $q$ 
is at most equal to the dimension of $\phi$, so 
$q$ is isomorphic to a form of type $h_{K}^{s} \oplus \phi$. 
By Lemma~\ref{lem:sum-spin} and Theorem~\ref{thm.spin} we get the inequalities
\[  \frac{3n}{7} \geq 3r \geq\ed_{n}(q) = \ed\rT_{n} \geq \ed \Spin_{n} - 1 
\, . \]
The resulting inequality fails for every even $n \ge 12$ because, for 
such $n$, $\ed \Spin_n \ge n/2$; see~\eqref{e.old-spin}. 

So we may assume that $7r \geq n$; then there is an isomorphism between 
the quadratic forms $\phi$ 
and a form of the type $h_{K}^{\oplus s}\oplus q$. 
By comparing dimensions we get the equality $7r = n + 2s$ when $r$ is even, and $7r+1 = n + 2s$ when $r$ is odd.
The essential dimension of the form 
$\phi$ as an element 
of $\rT_{7r}(K)$ or $\rT_{7r+1}(K)$ is at most $3r$, 
while by Lemma~\ref{lem:decrease-ed} we have that this essential 
dimension is at least $\ed_{n}(q) - s(s + 2n -1)/2$. From this, Lemma~\ref{lem.T-Spin} and Theorem~\ref{thm.spin} we have the chain of inequalities
   \begin{align*} 
   3r &\geq \ed_{n}(q) - \frac{s(s + 2n -1)}{2}\\
   &= \ed\rT_{n} - \frac{s(s + 2n -1)}{2}\\
   &\geq \ed\Spin_{n} - 1 - \frac{s(s + 2n -1)}{2}\\
   &\geq 2^{(n-2)/2} - \frac{n(n-1)}{2} - 1 - \frac{s(s + 2n -1)}{2}.
   \end{align*}

Now assume that $r$ is even. We substitute the value $s = (7r - n)/2$ in the resulting inequality, obtaining
   \[
   \frac{49r^{2} + (14n+10)r
      - 2^{(n+4)/2} - n^{2} + 2n - 8}{8}
   \geq 0.
   \]

We interpret this as a quadratic inequality in $r$. 
The constant term of the polynomial is negative
for all $n \geq 8$; hence if $r_{0}$ is the positive 
root, the equality is equivalent to $r \geq r_{0}$. 
By the quadratic formula 
   \begin{align*}
      r_0 &= \frac{\sqrt{49 \cdot 2^{(n+4)/2} + 168n - 367} -(7n + 5)}{49}\\
      &\geq  \frac{2^{(n+4)/4} -n-2}{7}\, .
   \end{align*}
This completes the proof of Theorem~\ref{thm.pfister} when $r$ is even. The calculations when $r$ is odd are analogous: using the substitution $s = (7r+1-n)/2$ we obtain the root
   \begin{align*}
      r_0 &= \frac{\sqrt{49 \cdot 2^{(n+4)/2} + 168n - 199} -(7n + 12)}{49}\\
      &\geq  \frac{2^{(n+4)/4} -n-2}{7}\, .
   \end{align*}
\end{proof}


\bibliographystyle{amsalpha}
\def\cprime{$'$}
\providecommand{\bysame}{\leavevmode\hbox to3em{\hrulefill}\thinspace}
\providecommand{\MR}{\relax\ifhmode\unskip\space\fi MR }
\providecommand{\MRhref}[2]{%
  \href{http://www.ams.org/mathscinet-getitem?mr=#1}{#2}
}
\providecommand{\href}[2]{#2}

\end{document}